\documentclass[12pt,a4paper,fleqn]{article}

\usepackage{lscape,exscale,amsthm}
\usepackage[intlimits]{amsmath}
\usepackage{rawfonts}
\usepackage{latexsym}
\usepackage[cp850]{inputenc}
\usepackage{epsfig}
\usepackage{bbm}
\usepackage{enumitem}
\usepackage{float}

\RequirePackage[OT1]{fontenc}
\RequirePackage{amsthm,amsmath}
\RequirePackage[numbers]{natbib}
\RequirePackage[colorlinks,citecolor=blue,urlcolor=blue]{hyperref}
\newcommand{\bol}[1]{\mbox{\boldmath$#1$}}

\newcommand{\bSigma}{\bol{\Sigma}}

\newcommand{\bx}{\mathbf{X}}

\newcommand{\by}{\mathbf{Y}}

\newcommand{\bH}{\mathbf{H}}

\newcommand{\bI}{\mathbf{I}}

\newcommand{\bA}{\bol{A}}

\newcommand{\bn}{\boldsymbol{\nu}}

\newcommand{\bxi}{\boldsymbol{\xi}}

\newcommand{\bS}{\mathbf{S}}
\newcommand{\sx}{\bar{\mathbf{x}}}
\newcommand{\Tr}{\text{tr}}

\newcommand{\bTheta}{\mathbf{\Theta}}
\newcommand{\ta}{\alpha}
\newcommand{\tb}{\beta}

\usepackage{stackrel}

\numberwithin{equation}{section}
\theoremstyle{plain}
\newtheorem{theorem}{Theorem}[section]

\newtheorem{proposition}{Proposition}[section]
\newtheorem{lemma}{Lemma}[section]
\newtheorem{corollary}{Corollary}[section]

\begin{document}

\begin{center}
\vspace*{2cm} \noindent {\bf \large Optimal Linear Shrinkage Estimator for Large Dimensional Precision Matrix}\\
\vspace{1cm} \noindent {\sc  Taras Bodnar$^{a}$, Arjun K. Gupta$^{b}$\footnote{Corresponding author. E-mail address: gupta@bgsu.edu. The first author is partly supported by the German Science Foundation (DFG) via the Research Unit 1735 ''Structural Inference in Statistics: Adaptation and Efficiency''.} and Nestor Parolya$^{c}$}\\
\vspace{1cm}
{\it \footnotesize  $^a$
Department of Mathematics, Humboldt-University of Berlin, D-10099 Berlin, Germany}\\
{\it \footnotesize  $^b$
Department of Mathematics and Statistics,  Bowling Green State University, Bowling Green, OH 43403, USA} \\
{\it \footnotesize  $^c$
Department of Statistics and Econometrics, Ruhr-University Bochum, 44801 Bochum, Germany} \\
\end{center}

\begin{abstract}
In this work we construct an optimal shrinkage estimator for the precision matrix in high dimensions. We consider the general asymptotics when the number of variables $p\rightarrow\infty$ and the sample size $n\rightarrow\infty$ so that $p/n\rightarrow c\in (0, +\infty)$. The precision matrix is estimated directly, without inverting the corresponding estimator for the covariance matrix. The recent results from the random matrix theory allow us to find the asymptotic deterministic equivalents of the optimal shrinkage intensities and estimate them consistently. The resulting distribution-free estimator has almost surely the minimum Frobenius loss.  Additionally, we prove that the Frobenius norms of the inverse and of the pseudo-inverse sample covariance matrices tend almost surely to deterministic quantities and estimate them consistently. At the end, a simulation is provided where the suggested estimator is compared with the estimators for the precision matrix proposed in the literature. The optimal shrinkage estimator shows significant improvement and robustness even for non-normally distributed data.
\end{abstract}

\vspace{0.7cm}


\section{Introduction}
 The estimation of the covariance matrix, as well as its inverse (the precision matrix), plays an important role in many disciplines from finance and genetics to wireless communications and engineering. In fact, having a suitable estimator for the precision matrix we are able to construct a good estimator for different types of optimal portfolios (see, Markowitz (1952), Elton et al. (2009)). Similarly, in the array processing, the beamformer or the so-called minimum variance distortionless response spatial filter is defined in terms of the precision matrix (see, e.g., Van Trees (2002)). In practice, however, the true precision matrix is unknown and a feasible estimator, constructed from data, must be used.

 If the number of variables $p$ is much smaller than the sample size $n$ we can use the sample estimator which is biased but a consistent estimator for the precision matrix (see, e.g., Bai and Shi (2011)). This case is known in the multivariate statistics as the "standard asymptotics" (see, Le Cam and Yang (2000)). There are many findings on the estimation of the precision matrix when a particular distribution assumption is imposed. For example, the estimation of the precision matrix under the multivariate normal distribution was considered by Krishnamoorthy and Gupta (1989), Gupta and Ofori-Nyarko (1994, 1995a, 1995b), Kubokawa (2005) and Tsukuma and Konno (2006). The results in the case of multivariate Pearson type II distribution as well as the multivariate elliptically contoured stable distribution are obtained by Sarr and Gupta (2009) as well as by Bodnar and Gupta (2011) and Gupta et al. (2013), respectively.

 Unfortunately, in practice $p$ is often comparable in size to $n$ or even is greater than $n$, i.e. we are in the situation when both the sample size $n$ and the dimension $p$ tend to infinity but their ratio keeps (tends to) a positive constant. This case often arises in finance when the number of assets is comparable or even greater than the number of observations for each asset. Similarly, in genetics, the data set can be huge comparable to the number of patients. Both examples illustrate the importance of the results obtained for $p,n\rightarrow\infty$.

 We deal with this type of asymptotics, called the "large dimensional asymptotics" and also known as the "Kolmogorov asymptotics", in the present paper. More precisely, it is assumed that the dimension $p\equiv p(n)$ is a function of the sample size $n$ and $p/n\rightarrow c\in(0, +\infty)$ as $n\rightarrow\infty$. This general type of asymptotics was intensively studied by several authors (see, Girko (1990, 1995), B\"{u}hlmann and van Geer (2011) etc.). In this asymptotics the usual estimators for the precision matrix perform poorly and are not consistent anymore. There are some techniques which can be used to handle the problem. Assuming that the covariance (precision) matrix has a sparse structure, significant improvements have already been achieved (see, Cai et al. (2011), Cai and Shen (2011), Cai and Zhou (2012)). For the low-rank covariance matrices see the work of Rohde and Tsybakov (2011). An interesting nonparanormal graphic model was recently proposed by Xue and Zou (2012). Also, in order to estimate the large dimensional covariance matrix the method of block thresholding can be applied (see, Cai and Yuan (2012)). If the covariance matrix has a factor structure then the progress has been made by Fan et al. (2008).

 However, if neither the assumption about the structure of covariance (precision) matrix nor about a particular distribution is imposed, only a few results are known in the literature which are based on the shrinkage estimators in high-dimensional setting (cf. Ledoit and Wolf (2004), Ledoit and Wolf (2012), Bodnar et al. (2013)). The shrinkage estimator was first developed by Stein (1956) and forms a linear combination of the sample estimator and some target. The corresponding shrinkage coefficients are often called shrinkage intensities. Ledoit and Wolf (2004) proposed to shrink the sample covariance matrix to the identity matrix and showed that the resulting estimator is well-behaved in large dimensions and it is optimal in the sense of quadratic mean. This estimator is called the linear shrinkage estimator because it shrinks the eigenvalues of the sample covariance matrix linearly. Recently, Bodnar et al. (2013) proposed a generalization of the linear shrinkage estimator, where the shrinkage target was chosen to be an arbitrary nonrandom matrix\footnote{Of course, as the covariance matrix is assumed to be positive definite and symmetric, the target matrix must also possess this properties.} and they showed the almost sure convergence of the derived estimator to its oracle. Some new results are obtained by Ledoit and Wolf (2012), who considered the so-called nonlinear shrinkage estimator and derived it using the theory of the random matrices.

The aim of our paper is to construct a feasible estimator for the precision matrix using the linear shrinkage technique and the random matrix theory. In contrast to well-known procedures, we shrink the inverse of the sample covariance matrix itself instead of shrinking the sample covariance matrix and then inverting it. In the case when $c>1$ the pseudo inverse of the sample covariance matrix is taken. The recent results from the random matrix theory allow us to find the asymptotics of the optimal shrinkage intensities and estimate them consistently.

The random matrix theory is a very fast growing branch of the probability theory with many applications in statistics. It studies the asymptotic behavior of the eigenvalues of the different random matrices under general asymptotics (see, e.g., Anderson et al. (2010), Bai and Silverstein (2010)). The asymptotic behavior of the functionals of the sample covariance matrices was studied by Ma$\breve{\text{c}}$enko and Pastur (1967), Yin (1986), Girko and Gupta (1994, 1996a, 1996b), Silverstein (1995), Bai et al. (2007), Bai and Silverstein (2010), Rubio and Mestre (2011) etc.

We extend these results in the present paper by establishing the almost sure convergence of the optimal shrinkage intensities and the Frobenius norm of the inverse sample covariance matrix.
Moreover, we construct a general linear shrinkage estimator for the precision matrix which has \textit{almost surely} the smallest Frobenius loss when both the dimension $p$ and the sample size $n$ increase together and $p/n\rightarrow c\in(0, +\infty)$ as $n\rightarrow\infty$.

The rest of the paper is organized as follows. In Section 2 we present some preliminary results from the random matrix theory and formulate the assumptions used throughout the paper. In Section 3 we construct the \textit{oracle} linear shrinkage estimator for the precision matrix and verify the main asymptotic results about the shrinkage intensities and the Frobenius norm of the inverse and pseudo-inverse sample covariance matrices. Section 4 is dedicated to the \textit{bona fide} linear shrinkage estimator for the precision matrix while Section 5 contains the results of the simulation study. Here, the performance of the derived estimator is compared with other known estimators for the large dimensional precision matrices. Section 6 includes the summary, while the proofs of the theorems are presented in the appendix (Section 7).

\section{Assumptions and notations}
The "large dimensional asymptotics" or "Kolmogorov asymptotics" include  $\dfrac{p}{n}\rightarrow c\in (0, +\infty)$ as both the number of variables $p\equiv p(n)$ and the sample size $n$ tend to infinity. In this case the traditional sample estimator performs poorly or very poorly and tends to over/underestimate the population inverse covariance matrix. The inverse of the sample covariance matrix $\bS^{-1}_n$ is biased, inconsistent for $\dfrac{p}{n}\rightarrow c>0$ as $n\rightarrow\infty$  and it does not exist for $c>1$. For example, under the normality assumption $\bS^{-1}_n$ has an inverse Wishart distribution if $c<1$, (cf. Gupta and Nagar (2000))
 \begin{equation*}
  E(\bS^{-1}_n)=\dfrac{n}{n-p-2}\bSigma^{-1}_n\,
 \end{equation*}
In particular, for $p=n/2+2$ we have that $c=1/2$ and $E(\bS^{-1}_n)=2\bSigma^{-1}_n$. In general, as $c$ increases  the sample estimator of the precision matrix becomes worse.

We use the following notations in the paper:
\begin{itemize}
\item $\bSigma_n$ stands for the covariance matrix, $\bS_n$ denotes the corresponding sample covariance matrix.\footnote{Since the dimension $p\equiv p(n)$ is a function of the sample size $n$, the covariance matrix $\bSigma_n$ also depends on $n$ via $p(n)$. That is why we make use of the subscript $n$ for all of the considered objects in order to emphasize this fact and to simplify the notation in the paper.} The population covariance matrix $\bSigma_n$ is a nonrandom $p$-dimensional positive definite matrix.
\item $||\bA||^2_F=\text{tr}(\bA\bA^\prime)$ denotes the Frobenius norm of a square matrix $\bA$, $||\bA||_{tr}=\text{tr}\left[\left(\bA\bA^\prime\right)^{1/2}\right]$ stands for its trace norm, while $||\bA||_{2}$  is the spectral norm.
\item The pairs $(\tau_i,\bn_i)$ for $i=1,\ldots,p$ denote the collection of eigenvalues and the corresponding orthonormal eigenvectors of the covariance matrix $\bSigma_n$.
\item $H_n(t)$ is the empirical distribution function (e.d.f.) of the eigenvalues of $\bSigma_n$, i.e.
\begin{equation}\label{cov_edf}
H_n(t)=\dfrac{1}{p}\sum\limits_{i=1}^{p}\mathbbm{1}_{\{\tau_i\leq t\}}\,
\end{equation}
 where $\mathbbm{1}_{\{\cdot\}}$ is the indicator function.
\item Let $\bx_n$ be a $p\times n$ matrix which consists of independent and identically distributed (i.i.d.) real random variables with zero mean and unit variance. The observation matrix is defined as
\begin{equation}\label{obs}
 \by_n=\bSigma_n^{\frac{1}{2}}\bx_n.
 \end{equation}
Only the matrix $\by_n$ is observable. We know neither $\bx_{n}$ nor $\bSigma_n$ itself.
\item The pairs $(\lambda_i,\mathbf{u}_i)$ for $i=1,\ldots,p$ are the eigenvalues and the corresponding orthonormal eigenvectors of the sample covariance matrix\footnote{The sample mean vector $\sx$ was omitted because the $1$-rank matrix $\sx\sx^{\prime}$ does not influence the asymptotic behavior of the spectrum of sample covariance matrix (see, Bai and Silverstein (2010), Theorem A.44).}
 \begin{equation}\label{samplecov}
 \bS_n=\dfrac{1}{n}\by_n\by_n^{\prime}=\dfrac{1}{n}\bSigma_n^{\frac{1}{2}}\bx_n\bx_n^{\prime}\bSigma_n^{\frac{1}{2}}\,.
 \end{equation}
\item Similarly, the (e.d.f.) of the eigenvalues of the sample covariance matrix $\bS_n$ is defined as
\begin{equation}\label{sample_edf}
F_n(\lambda)=\dfrac{1}{p}\sum\limits_{i=1}^{p}\mathbbm{1}_{\{\lambda_i\leq\lambda\}}~~\forall~\lambda\in\mathbbm{R}\,.
\end{equation}
\item In order to handle the case when $c>1$ we introduce the dual sample covariance matrix defined as
\begin{equation}\label{dualsamplecov}
 \bar{\bS}_n=\dfrac{1}{n}\by_n^{\prime}\by_n=\dfrac{1}{n}\bx_n^{\prime}\bSigma_n\bx_n\,
 \end{equation}
 with the corresponding (e.d.f) defined by
\begin{equation}\label{dualsample_edf}
\bar{F}_n(\lambda)=\dfrac{1}{n}\sum\limits_{i=1}^{n}\mathbbm{1}_{\{\lambda_i<\lambda\}}~~\forall~\lambda\in\mathbbm{R}\,.
\end{equation}
 Note that the matrix $\bar{\bS}_n$ has same nonzero eigenvalues as $\bS_n$, they differ only in $|p-n|$ zero eigenvalues.

\end{itemize}

The main assumptions, which we mention throughout the paper, are as follows
\begin{description}
\item[(A1)] We assume that $H_n(t)$ converges to a limit $H(t)$ at all points of continuity of $H$.

\item[(A2)] The elements of the matrix $\bx_n$ have uniformly bounded $4+\varepsilon,~\varepsilon>0$ moments.

\item[(A3)] For all $n$ large enough there exists the compact interval $[h_0, h_1]$ in $(0,+\infty)$ which contains the support of $H_n$.
\end{description}

All of these assumptions are quite general and are satisfied in many practical situations. The assumption (A1) is essential to prove the Mar$\breve{\text{c}}$henko-Pastur equation (see, e.g., Silverstein (1995)) which is  used for studying the asymptotic behavior of the spectrum of general random matrices (see, e.g., Bai and Silverstein (2010)). The fourth moment is needed for the proof of Theorem 3.2 and Theorem 3.3. The assumption (A3) ensures that both the matrix $\bSigma_n$ and its inverse $\bSigma^{-1}_n$ have uniformly bounded spectral norms at infinity. It means that $\bSigma_n$ has the uniformly bounded maximum eigenvalue and its minimum eigenvalue is greater than zero. Rubio et al. (2012) pointed out that (A2) and (A3) are only some technical conditions which can be further violated.

In order to investigate the (e.d.f) $F_n(\lambda)$ the Stieltjes transform is used. For nondecreasing function with bounded variation $G$ the Stieltjes transform is defined as
\begin{equation}\label{Strans}
\forall z\in\mathbbm{C}^+~~~m_G(z)=\int\limits_{-\infty}^{+\infty}\dfrac{1}{\lambda-z}dG(\lambda)\,.
\end{equation}
In our notation  $\mathbbm{C}^+=\{z\in\mathbbm{C}: \textbf{Im}(z)>0\}$ is the half-plane of complex numbers with strictly positive imaginary part and any complex number is defined as $z=\textbf{Re}(z)+i\textbf{Im}(z)$. More about the Stieltjes transform and its properties can be found in Silverstein (2009).

The Stieltjes transform of the sample (e.d.f.) $F_n(\lambda)$ for all $z\in\mathbbm{C}^+$ is given by
\begin{eqnarray}\label{trans_s}
m_{F_n}(z)&=&\dfrac{1}{p}\sum\limits_{i=1}^{p}\int\limits_{-\infty}^{+\infty}\dfrac{1}{\lambda-z}\delta(\lambda-\lambda_i)d\lambda=\dfrac{1}{p}\Tr\{(\bS_n-z\bI)^{-1}\}\,
\end{eqnarray}
where $\bI$ is a suitable identity matrix and $\delta(\cdot)$ is the Dirac delta function.

\section{Optimal linear shrinkage estimator for the precision matrix}

\subsection{Case $c<1$}
In this section we construct an optimal linear shrinkage estimator for the precision matrix under high-dimensional asymptotics. The estimator is an \textit{oracle} one, i.e., it depends on unknown quantities. The corresponding \textit{bona fide} estimator is given in Section 4.
   We use a procedure similar to Bodnar et al. (2013) where the optimal linear shrinkage estimator for the covariance matrix was constructed.
 The general linear shrinkage estimator of the precision matrix $\bSigma_n^{-1}$ for $c<1$ is given by
\begin{equation}\label{gse}
\widehat{\boldsymbol{\Pi}}_{GSE}=\alpha_n\bS_n^{-1}+\beta_n\boldsymbol{\Pi}_0~~\text{with}~ \sup\limits_{p}||\boldsymbol{\Pi}_0||_{tr}\leq M\,.
\end{equation}
Note that we need the condition $c<1$ to keep the sample covariance matrix $\bS_n$ invertible. The assumption that the target matrix $\boldsymbol{\Pi}_0$ has a uniformly bounded trace norm, i.e. there exists $M>0$ such that $\sup\limits_{p}||\boldsymbol{\Pi}_0||_{tr}\le M$, is rather general and it is actually needed to keep the coefficient $\beta_n$ bounded for large dimensions $p$. This condition can be replaced with an equivalent assumption on $\beta_n$. Note that the target matrix can also be random but independent of $\by_{n}$.\footnote{In practice, however, one has to be careful with the choice of the target matrix $\boldsymbol{\Pi}_0$. If it is close in some sense to $\bS^{-1}_n$, negative shrinkage intensities might occur.}

Our aim is now to find the optimal shrinkage intensities which minimize the Frobenius-norm loss for a given nonrandom target matrix $\boldsymbol{\Pi}_0$ expressed as
\begin{equation}\label{risk1}
L^2_F=||\widehat{\boldsymbol{\Pi}}_{GSE}-\bSigma^{-1}_n||_F^2=||\bSigma^{-1}_n||^2_F+||\widehat{\boldsymbol{\Pi}}_{GSE}||^2_F-2\text{tr}
\left(\widehat{\boldsymbol{\Pi}}_{GSE}\bSigma^{-1}_n\right)\,,
\end{equation}

As a result, using (\ref{gse}) the following optimization problem has to be solved
\begin{eqnarray}\label{optm1}
&&\ta_n^2||\bS^{-1}_n||^2_F+2\ta_n\tb_n\text{tr}(\bS^{-1}_n\boldsymbol{\Pi}_0)+\tb_n^2||\boldsymbol{\Pi}_0||^2_F
-2\ta_n\text{tr}(\bS^{-1}_n\bSigma^{-1}_n)-2\tb_n\text{tr}(\bSigma^{-1}_n\boldsymbol{\Pi}_0)\longrightarrow\text{min}\nonumber\\
&&~~\text{with respect to}~\ta_n~\text{and}~\tb_n\nonumber\,.
\end{eqnarray}

Next, taking the derivatives of $L^2_F$ with respect to $\ta_n$ and $\tb_n$ and setting them equal to zero we get
\begin{equation}\label{dera}
\dfrac{\partial L^2_F}{\partial\ta_n}=\ta_n||\bS^{-1}_n||^2_F+\tb_n\text{tr}(\bS^{-1}_n\boldsymbol{\Pi}_0)-\text{tr}(\bS^{-1}_n\bSigma^{-1}_n)=0\,,
\end{equation}
\begin{equation}\label{derb}
\dfrac{\partial L^2_F}{\partial\tb_n}=\ta_n\text{tr}(\bS^{-1}_n\boldsymbol{\Pi}_0)+\tb_n||\boldsymbol{\Pi}_0||^2_F-\text{tr}(\bSigma^{-1}_n\boldsymbol{\Pi}_0)=0\,.
\end{equation}
The Hessian of the $L^2_F$ has the form
\begin{equation}\label{hessian1}
 \bH=\left(
 \begin{array}{cc}
 ||\bS^{-1}_n||^2_F&\text{tr}(\bS^{-1}_n\boldsymbol{\Pi}_0)\\
 \text{tr}(\bS^{-1}_n\boldsymbol{\Pi}_0)&||\boldsymbol{\Pi}_0||^2_F\,
 \end{array}
 \right)
\end{equation}
which is always positive definite, since
\begin{align}\label{posdef}
 \text{det}(\bH)&=||\bS^{-1}_n||^2_F||\boldsymbol{\Pi}_0||^2_F-(\text{tr}(\bS^{-1}_n\boldsymbol{\Pi}_0))^2\\
 &\geq||\bS^{-1}_n||^2_F||\boldsymbol{\Pi}_0||^2_F-||\bS^{-1}_n||^2_{2}(\text{tr}(\boldsymbol{\Pi}_0))^2
 \stackrel{Jensen}{\geq}(||\bS^{-1}_n||^2_F-||\bS^{-1}_n||^2_{2})||\boldsymbol{\Pi}_0||^2_F>0\nonumber\,,
\end{align}
where the last inequality in (\ref{posdef}) is well-known (see, e.g., Horn and Johnson (1985)).

 Thus, the optimal $\alpha^*_n$ and $\beta^*_n$ are given by

 \begin{equation}\label{talfa} \alpha_n^*=\dfrac{\text{tr}(\bS^{-1}_n\bSigma^{-1}_n)||\boldsymbol{\Pi}_0||^2_F-\text{tr}(\bSigma^{-1}_n\boldsymbol{\Pi}_0)\text{tr}(\bS^{-1}_n\boldsymbol{\Pi}_0)}{||\bS^{-1}_n||^2_F||\boldsymbol{\Pi}_0||^2_F-\bigl(\text{tr}(\bS^{-1}_n\boldsymbol{\Pi}_0)\bigr)^2}\,,
 \end{equation}
 \begin{equation}\label{tbeta} \beta_n^*=\dfrac{\text{tr}(\bSigma^{-1}_n\boldsymbol{\Pi}_0)||\bS^{-1}_n||^2_F-\text{tr}(\bS^{-1}_n\bSigma^{-1}_n)\text{tr}(\bS^{-1}_n\boldsymbol{\Pi}_0)}{||\bS^{-1}_n||^2_F||\boldsymbol{\Pi}_0||^2_F-\bigl(\text{tr}(\bS^{-1}_n\boldsymbol{\Pi}_0)\bigr)^2}\,.
 \end{equation}

Now, we formulate our first main result in Theorem 3.1 which states that the normalized Frobenius norm of the inverse sample covariance matrix $1/p||\bS^{-1}_n||^2_F$ tends almost surely to a nonrandom quantity.
\begin{theorem}
Assume that (A1) and (A3) hold and $\dfrac{p}{n}\rightarrow c\in(0, 1)$ for $n\rightarrow\infty$. Then the normalized Frobenius norm of the inverse sample covariance matrix $\psi_n=\dfrac{1}{p}||\bS^{-1}_n||^2_F$ almost surely tends to a nonrandom $\psi$ which is given by
\begin{equation}\label{psi} \psi=\dfrac{1}{(1-c)^2}\int\limits_{-\infty}^{+\infty}\dfrac{dH(\tau)}{\tau^2}+\dfrac{c}{(1-c)^3}\left(\int\limits_{-\infty}^{+\infty}\dfrac{dH(\tau)}{\tau}\right)^2\,.
\end{equation}
\end{theorem}
The proof is given in the Appendix. Theorem 3.1 presents an important result which indicates that the Frobenius norm of the inverse sample covariance matrix is asymptotically nonrandom as well as it depends on $H$ and concentration ratio $c$ only. Moreover, Theorem 3.1 gives us an intuitive hint how to find the asymptotic equivalent representation of $||\bS^{-1}_n||^2_F$. The corresponding result is presented in Theorem 3.2.
\begin{theorem}
Let the assumptions (A1)-(A3) hold and $\dfrac{p}{n}\rightarrow c\in(0, 1)$. Then as $n\rightarrow\infty$,
\begin{equation}\label{phi2}
\dfrac{1}{p}\Biggl|||\bS^{-1}_n||^2_F-\left(\dfrac{1}{(1-c)^2}||\bSigma^{-1}_n||^2_F+\dfrac{c}{p(1-c)^3}||\bSigma^{-1}_n||_{tr}^2\right)\Biggr|\underset{\text{a.s.}}{\longrightarrow}0\,.
\end{equation}

Additionally, for the quantity $\text{tr}(\bS^{-1}_n\mathbf{\Theta})$ with a symmetric positive definite matrix $\mathbf{\Theta}$ which has uniformly bounded trace norm as $n\rightarrow\infty$,
\begin{equation}\label{trst}
\Biggl|\text{tr}(\bS^{-1}_n\mathbf{\Theta})-\dfrac{1}{1-c}\text{tr}(\bSigma^{-1}_n\mathbf{\Theta})\Biggr|\underset{\text{a.s.}}{\longrightarrow}0~~\text{for}~~\dfrac{p}{n}\rightarrow c\in(0, 1)\,.
\end{equation}
\end{theorem}
The theorem is proved in the Appendix. Theorem 3.2 provides us the asymptotic behavior of the Frobenius norm of the inverse sample covariance matrix and of the functional $\text{tr}(\bS^{-1}_n\mathbf{\Theta})$. It shows that the consistent estimator for the Frobenius norm of the precision matrix under the general asymptotics is not equal to its sample counterpart. Using Theorem 3.2 we can easily determine the asymptotic bias of the sample estimator which consists of the two types of biases. The multiplicative bias is violated by multiplying $||\bS^{-1}_n||^2_F$ by $(1-c)^2$. After that, the additive bias is dealt by subtracting $\dfrac{c}{p(1-c)}||\bSigma^{-1}_n||_{tr}^2$.  The sample estimator of the functional $\text{tr}(\bS^{-1}_n\mathbf{\Theta})$ is also not a consistent estimator for $\text{tr}(\bSigma^{-1}_n\mathbf{\Theta})$. The consistent estimator is obtained by multiplying $\text{tr}(\bS^{-1}_n\mathbf{\Theta})$ by the constant $(1-c)$.

Results similar to those given in Theorem 3.1 and Theorem 3.2 are also available for the estimation of the population covariance matrix (cf. Bodnar et al. (2013)). However, in the case of the covariance matrix, the sample estimator for the Frobenius norm possesses only the additive bias $\dfrac{c}{p}\text{tr}(||\bSigma_n^{-1}||_{tr})$, while $\text{tr}(\bS_n\bTheta)$ is a consistent estimator for $\text{tr}(\bSigma_n\bTheta)$.

Next, we show that the optimal shrinkage intensities $\alpha^*_n$ and $\beta^*_n$ are almost surely asymptotic equivalent to nonrandom quantities $\alpha^*$ and $\beta^*$ under the large-dimensional asymptotics $\dfrac{p}{n}\rightarrow c\in(0, 1)$.

\begin{corollary}
Assume that (A1)-(A3) hold and $\dfrac{p}{n}\rightarrow c\in(0, 1)$ for $n\rightarrow\infty$. Then for the optimal shrinkage intensities $\ta^*_n$ and $\tb^*_n$
\begin{equation}\label{ainfty}
\left|\alpha^*_n-\alpha^*\right|\longrightarrow0~~\text{a. s.}\,,
\end{equation}
where
\begin{equation}\label{a_as}
\alpha^*=(1-c)\dfrac{||\bSigma_n^{-1}||^2_F||\boldsymbol{\Pi}_0||^2_F-\left(\text{tr}(\bSigma^{-1}_n\boldsymbol{\Pi}_0)\right)^2}{\left(||\bSigma_n^{-1}||^2_F+\dfrac{c}{p(1-c)}||\bSigma_n^{-1}||_{tr}^2\right)||\boldsymbol{\Pi}_0||^2_F-\left(\text{tr}(\bSigma^{-1}_n\boldsymbol{\Pi}_0)\right)^2}
\end{equation}
and
\begin{equation}\label{binfty}
\left|\beta^*_n-\beta^*\right|\longrightarrow0~~\text{a. s.}\,,
\end{equation}
with
\begin{equation}\label{b_as}
\beta^*=\dfrac{\text{tr}(\bSigma^{-1}_n\boldsymbol{\Pi}_0)}{||\boldsymbol{\Pi}_0||^2_F}\left(1-\dfrac{\alpha^*}{1-c}\right)
\end{equation}
\end{corollary}

Note that both the asymptotic optimal intensities $\ta^*$ and $\tb^*$ are always positive as well as $\ta^*\in(0,1-c)$ due to inequality (\ref{posdef}) and $c\in(0,1)$. Using these results we are immediately able to estimate $\ta^*$, $\tb^*$ consistently which is shown in Section 4.

\subsection{Case $c>1$}

In this subsection we deal with the problem of the estimation of the precision matrix when the dimension $p$ is greater than the sample size $n$, i.e., $c>1$. This case is very difficult to handle because of the loss of information as $c$ becomes greater than one. Moreover, the sample covariance matrix $\bS_n$ is not invertible and thus the estimator $\bS_n^{-1}$ must be replaced by a suitable one. This is usually done by using the generalized inverse matrix $\bS_n^{+}$ instead of $\bS_n^{-1}$. In this case the general shrinkage estimator has the form
\begin{equation}\label{gse1}
\widehat{\boldsymbol{\Pi}}_{GSE}=\tilde{\alpha}_n\bS_n^{+}+\tilde{\beta}_n\boldsymbol{\Pi}_0~~\text{with}~ \sup\limits_{p}||\boldsymbol{\Pi}_0||_{tr}\leq M\,.
\end{equation}
The optimal shrinkage intensities $\tilde{\alpha}^{*}_n$ and $\tilde{\beta}^{*}_n$ are determined following the procedure of Section 3.1. They are given by
\begin{equation}\label{talfa1} \tilde{\alpha}_n^*=\dfrac{\text{tr}(\bS^{+}_n\bSigma^{-1}_n)||\boldsymbol{\Pi}_0||^2_F-\text{tr}(\bSigma^{-1}_n\boldsymbol{\Pi}_0)\text{tr}(\bS^{+}_n\boldsymbol{\Pi}_0)}{||\bS^{+}_n||^2_F||\boldsymbol{\Pi}_0||^2_F-\bigl(\text{tr}(\bS^{+}_n\boldsymbol{\Pi}_0)\bigr)^2}\,,
 \end{equation}
 \begin{equation}\label{tbeta1} \tilde{\beta}_n^*=\dfrac{\text{tr}(\bSigma^{-1}_n\boldsymbol{\Pi}_0)||\bS^{+}_n||^2_F-\text{tr}(\bS^{+}_n\bSigma^{-1}_n)\text{tr}(\bS^{+}_n\boldsymbol{\Pi}_0)}{||\bS^{+}_n||^2_F||\boldsymbol{\Pi}_0||^2_F-\bigl(\text{tr}(\bS^{+}_n\boldsymbol{\Pi}_0)\bigr)^2}\,.
 \end{equation}

In Theorem 3.3 we derive the asymptotic properties of two quantities used in (\ref{talfa1}) and (\ref{tbeta1} ), namely $\text{tr}(\bTheta\bS^{+}_n)$ and $||\bS^{+}_n||^2_F$.

\begin{theorem}
Let the assumptions (A1)-(A3) hold and $\dfrac{p}{n}\rightarrow c\in(1, +\infty)$. Then as $n\rightarrow\infty$,
\begin{equation}\label{phi21}
\Biggl|\dfrac{1}{p}||\bS^{+}_n||^2_F-c^{-1}x^\prime(0)\Biggr|\underset{\text{a.s.}}{\longrightarrow}0,~~\text{where}~x^\prime(0)=\dfrac{1}{\dfrac{1}{x^2(0)}-\dfrac{c}{p}\text{tr}\left[\left(\bSigma^{-1}_n+x(0)\bI\right)^{-2}\right]}\,
\end{equation}
 and $x(0)$ is the unique solution of the equation
\begin{equation}\label{x01}
\dfrac{1}{x(0)}=\dfrac{c}{p}\text{tr}\left[\left(\bSigma^{-1}_n+x(0)\bI\right)^{-1}\right]
\end{equation}

Additionally, for the quantity $\text{tr}(\bTheta\bS^{+}_n)$ with a symmetric positive definite matrix $\mathbf{\Theta}$ which has uniformly bounded spectral norm, as $n\rightarrow\infty$,
\begin{equation}\label{trst1}
\Biggl|\dfrac{1}{p}\text{tr}(\bTheta\bS^{+}_n)-c^{-1}y(\bTheta)\Biggr|\underset{\text{a.s.}}{\longrightarrow}0~~\text{for}~~\dfrac{p}{n}\rightarrow c\in(1, +\infty)\,,
\end{equation}
where $y(\bTheta)$ is the solution of
\begin{equation}\label{y0}
\dfrac{1}{y(\bTheta)}=\dfrac{c}{p}\text{tr}\left[\left(\bSigma_n^{-1/2}
\bTheta\bSigma_n^{-1/2}+y(\bTheta)\bI\right)^{-1}\right]
\end{equation}
\end{theorem}
The proof of Theorem 3.3 is given in the Appendix. The results of Theorem 3.3 show that using the generalized inverse technique it is not clear how to estimate the functionals of $\bSigma_n^{-1}$ consistently. The asymptotic values obtained in Theorem 3.3 are far away from the desired ones. In order to correct these biases, we need to solve the non-linear equations (\ref{x01}) and (\ref{y0}), respectively, which appears to be a difficult task. Finally, we notice, that the quantities $x(0)$ and $x^\prime(0)$, however, can be estimated consistently using Theorem 3.3.

In an important special case when the matrix $\bTheta=\boldsymbol{\xi}\boldsymbol{\eta}^\prime$ for some $\bxi$ and $\boldsymbol{\eta}$ with bounded Euclidean norms we get the following result summarized in Proposition 3.1 which is proved in the Appendix.

\begin{proposition}
Under the assumptions of Theorem 3.3 and $\bTheta=\boldsymbol{\xi}\boldsymbol{\eta}^\prime$ for some $\bxi$ and $\boldsymbol{\eta}$ with bounded Euclidean norms it holds
\begin{equation}\label{special}
\dfrac{1}{p}\Biggl|\boldsymbol{\eta}^\prime\bS^{+}_n\boldsymbol{\xi}-\dfrac{c^{-1}}{c-1}\boldsymbol{\eta}^\prime\bSigma^{-1}_n\boldsymbol{\xi}\Biggr|\underset{\text{a.s.}}{\longrightarrow}0
~~\text{for}~~\dfrac{p}{n}\rightarrow c\in(1, +\infty)\,,
\end{equation}
\end{proposition}

It is remarkable to note that the results of Proposition 3.1 are very similar to those presented in Theorem 3.2 if $\bTheta=\boldsymbol{\xi}\boldsymbol{\xi}^\prime$. The only difference is the sign of the constant $(1-c)$.

Next we use the asymptotic results of Theorem 3.3 for finding the asymptotic equivalents to the optimal shrinkage intensities $\tilde{\alpha}^*_n$ and $\tilde{\beta}^{*}_n$ given in (\ref{talfa1}) and (\ref{tbeta1}), respectively.

\begin{corollary}
Assume that (A1)-(A3) hold and $\dfrac{p}{n}\rightarrow c\in(1, +\infty)$ for $n\rightarrow\infty$. Then for the optimal shrinkage intensities $\ta^*_n$ and $\tb^*_n$ from (\ref{talfa1}) and (\ref{tbeta1}) holds
\begin{equation}\label{ainfty11}
\left|\alpha^*_n-\alpha^*\right|\longrightarrow0~~\text{a. s.}\,,
\end{equation}
where
\begin{equation}\label{a_as1}
\alpha^*=\dfrac{y(\bSigma^{-1}_n)||\boldsymbol{\Pi}_0||^2_F-y(\boldsymbol{\Pi}_0)\text{tr}(\bSigma^{-1}_n\boldsymbol{\Pi}_0)}{x^\prime(0)||\boldsymbol{\Pi}_0||^2_F-c^{-1}y^2(\boldsymbol{\Pi}_0)}
\end{equation}
and
\begin{equation}\label{binfty11}
\left|\beta^*_n-\beta^*\right|\longrightarrow0~~\text{a. s.}\,,
\end{equation}
with
\begin{equation}\label{b_as1}
\beta^*=\dfrac{\text{tr}(\bSigma^{-1}_n\boldsymbol{\Pi}_0)x^\prime(0)-y(\bSigma^{-1}_n)y(\boldsymbol{\Pi}_0)}{cx^\prime(0)||\boldsymbol{\Pi}_0||^2_F-y^{2}(\boldsymbol{\Pi}_0)}
\end{equation}
\end{corollary}

Even if the target matrix $\boldsymbol{\Pi}_0$ is chosen as a one-rank matrix, i.e. $\boldsymbol{\Pi}_0=\boldsymbol{\xi}\boldsymbol{\eta}^\prime$, we are not able to provide consistent estimates for $\alpha^*$ and $\beta^*$ without an additional assumption imposed on $\bSigma_n$. One of possible assumptions for which $\alpha^*$ and $\beta^*$ are consistently estimable is $\bSigma_n=\sigma\bI$ as illustrated in Corollary 3.4 below. If $\bSigma_n=\sigma\bI$, then for $\dfrac{1}{p}||\bS^{+}_n||^2_F$ and $\dfrac{1}{p}\text{tr}(\bS^{+}_n)$ we get

\begin{corollary} Under the assumptions of Theorem 3.3 assume additionally that $\bSigma_n=\sigma\bI$. Then it holds as $n\rightarrow\infty$
\begin{equation}\label{spec1}
\Biggl|\dfrac{1}{p}||\bS^{+}_n||^2_F-\dfrac{\sigma^{-2}}{(c-1)^3}\Biggr|\underset{\text{a.s.}}{\longrightarrow}0\,,
\end{equation}
Additionally, for the quantity $\text{tr}(\bS^{+}_n)$ as $n\rightarrow\infty$ the norm
\begin{equation}\label{spec2}
\Biggl|\dfrac{1}{p}\text{tr}(\bS^{+}_n)-\dfrac{c^{-1}}{(c-1)}\sigma^{-1}\Biggr|\underset{\text{a.s.}}{\longrightarrow}0~~\text{for}~~\dfrac{p}{n}\rightarrow c\in(1, +\infty)\,.
\end{equation}
\end{corollary}

The proof of Corollary 3.3 is based on the fact that the equation (\ref{x01}) has the explicit solution $x(0)=\dfrac{\sigma^{-1}}{c-1}$ if $\bSigma_n=\sigma\bI$. The rest calculations are only technical ones. It is interesting to note that the result of Corollary 3.3 coincides with the corresponding one of Theorem 3.2 for $c<1$ if $\bSigma_n=\sigma\bI$.

Next we apply Corollary 3.3 with $\bSigma_n=\sigma\bI$ and $\boldsymbol{\Pi}_0=1/p\bI$ to construct the asymptotic equivalents to the optimal shrinkage intensities $\tilde{\alpha}^*_n$ and $\tilde{\beta}^{*}_n$  given in (\ref{talfa1}) and (\ref{tbeta1}), respectively.
\begin{corollary}
Assume that (A1)-(A3) hold, $\bSigma_n=\sigma\bI$, $\boldsymbol{\Pi}_0=1/p\bI$ and $\dfrac{p}{n}\rightarrow c\in(1, +\infty)$ for $n\rightarrow\infty$. Then for the optimal shrinkage intensities $\ta^*_n$ and $\tb^*_n$ from (\ref{talfa1}) and (\ref{tbeta1}) holds
\begin{equation}\label{alandbe}
\tilde{\alpha}^*_n\longrightarrow0~~\text{a. s.}~~\text{and}~~1/p\tilde{\beta}^*_n\longrightarrow\sigma^{-1}~~\text{a. s.}\,,
\end{equation}
\end{corollary}
Corollary 3.4 implies that the oracle optimal shrinkage estimator for the precision matrix in the case $c>1$ and $\bSigma_n=\sigma\bI$ is equal to
\begin{equation}\label{gse1000}
\widehat{\boldsymbol{\Pi}}_{GSE}=\bSigma^{-1}_n=\sigma^{-1}\bI\,.
\end{equation}
 The quantity $\sigma^{-1}=\dfrac{1}{p}\text{tr}(\bSigma^{-1}_n)$ can be easily estimated using the result of Corollary 3.3. Namely, the consistent estimator of $\sigma^{-1}$ is given by
\begin{equation}\label{sigma-1}
\hat{\sigma}^{-1}=p/n\dfrac{p/n-1}{p}\text{tr}(\bS^+_n)\,.
\end{equation}

However, in the general case when $\bSigma_n$ and $\boldsymbol{\Pi}_0$ are arbitrary, the results of Corollary 3.3 and Corollary 3.4 do not hold anymore. For this reason, we consider the oracle estimator given by (\ref{gse1}) in the simulation study of Section 5.

\section{Estimation of unknown parameters}
In this section we present consistent estimators for the asymptotic optimal shrinkage intensities derived in Section 3. The results of Theorem 3.2 allow us to estimate consistently the functionals of type $\text{tr}(\bSigma^{-1}_n\bTheta)$ and the Frobenius norm of the precision matrix. The consistent estimator for the functional $\theta_n(\bTheta)=\text{tr}(\bSigma^{-1}_n\bTheta)$ is given by
\begin{equation}\label{funcest}
\hat{\theta}_n(\bTheta)=(1-p/n)\text{tr}(\bS^{-1}_n\bTheta)\,
\end{equation}
which is a generalization of the so-called $G3$-estimator obtained by Girko (1995). In particular, in the case when $\bTheta=\boldsymbol{\xi}\boldsymbol{\eta}^\prime$ for some vectors $\boldsymbol{\xi}$ and $\boldsymbol{\eta}$ with bounded Euclidean norm, Girko (1995) showed that the corresponding estimator $\hat{\theta}_n(\boldsymbol{\xi}\boldsymbol{\eta}^\prime)$ tends to $\theta_n(\boldsymbol{\xi}\boldsymbol{\eta}^\prime)$ in probability. In contrast, Theorem 3.2 ensures the consistency of $\hat{\theta}_n(\bTheta)$ for more general forms of $\bTheta$ which should not be of rank 1.

Again, using (\ref{funcest}) and Theorem 3.2 we construct a consistent estimator for $\rho_n=1/p||\bSigma^{-1}_n||^2_F$ which is given by
\begin{equation}\label{Fnormest}
\hat{\rho}_n=\dfrac{(1-p/n)^2}{p}||\bS^{-1}_n||^2_F-\dfrac{1-p/n}{pn}||\bS^{-1}_n||_{tr}^2\,.
\end{equation}
Note that the result (\ref{Fnormest}) is entirely new and it was not mentioned in the literature up to now. Moreover, it is noted that for the derivation of (\ref{Fnormest}) we do not need the existence of the $4$th moment (see, the assumption (A2) in Section 2).

Using both the estimators (\ref{funcest}) and (\ref{Fnormest}), we are able now to construct the optimal linear shrinkage estimator (OLSE) for the precision matrix which is given by
\begin{equation}\label{olse}
\widehat{\boldsymbol{\Pi}}_{OLSE}=\widehat{\alpha}^*_n\bS_n^{-1}+\widehat{\beta}^*_n\boldsymbol{\Pi}_0~~\text{with}~ \sup\limits_{p}||\boldsymbol{\Pi}_0||_{tr}\leq M\,,
\end{equation}
where
\begin{align}\label{ea}
\widehat{\alpha}^*_n&=(1-p/n)\dfrac{p\hat{\rho}_n||\boldsymbol{\Pi}_0||^2_F-\hat{\theta}^2_n(\boldsymbol{\Pi}_0)}{\left(p\hat{\rho}_n+\dfrac{p/n}{p(1-p/n)}\hat{\theta}^2_n(\bI)\right)||\boldsymbol{\Pi}_0||^2_F-\hat{\theta}^2_n(\boldsymbol{\Pi}_0)}\nonumber\\
&=(1-p/n)\left(1-\dfrac{\dfrac{p/n}{p(1-p/n)}\hat{\theta}^2_n(\bI)||\boldsymbol{\Pi}_0||^2_F}{\left(p\hat{\rho}_n+\dfrac{p/n}{p(1-p/n)}\hat{\theta}^2_n(\bI)\right)||\boldsymbol{\Pi}_0||^2_F-\hat{\theta}^2_n(\boldsymbol{\Pi}_0)}\right)\nonumber\\
&=1-p/n-\dfrac{\dfrac{1}{n}||\bS^{-1}_n||_{tr}^2||\boldsymbol{\Pi}_0||^2_F}{||\bS_n^{-1}||^2_F||\boldsymbol{\Pi}_0||^2_F-\left(\text{tr}(\bS^{-1}_n\boldsymbol{\Pi}_0)\right)^2}\,.
\end{align}
and
\begin{equation}\label{eb}
\widehat{\beta}^*_n=\dfrac{\hat{\theta}_n(\boldsymbol{\Pi}_0)}{||\boldsymbol{\Pi}_0||^2_F}\left(1-\dfrac{\widehat{\alpha}_n^*}{1-p/n}\right)=\dfrac{\text{tr}(\bS^{-1}_n\boldsymbol{\Pi}_0)}{||\boldsymbol{\Pi}_0||^2_F}\left(1-p/n-\widehat{\alpha}_n^*\right)\,.
\end{equation}

The bona fide OLSE estimator (\ref{olse}) is optimal in the sense that it minimizes the Frobenius loss. It means that the estimators $\widehat{\alpha}_n^*$ and $\widehat{\beta}^*_n$ tend almost surely to their oracle asymptotic values (\ref{a_as}) and (\ref{b_as}) as $n\rightarrow\infty$, respectively. According to Corollary 3.1 the oracle optimal intensities $\ta^*_n$ and $\tb^*_n$ given in (\ref{talfa}) and (\ref{tbeta}) behave similarly. It is a remarkable property of the OLSE estimator which indicates that the bona fide OLSE estimator tends almost surely to its oracle one. Moreover, using the inequality (\ref{posdef}) it can be easily verified that the estimated optimal shrinkage intensities $\widehat{\alpha}^*_n$ and $\widehat{\beta}^*_n$ are almost surely positive and $\widehat{\alpha}^*_n$ has the support $(0,1-p/n)$. Only in the case when $p/n\rightarrow c=0$ as $n\rightarrow\infty$ the shrinkage intensities satisfy $\widehat{\alpha}^*_n \rightarrow 1$ and $\widehat{\beta}^*_n \rightarrow 0$. In this case the OLSE estimator coincides with the sample estimator which is consistent for the standard asymptotics.

If we compare the estimates of the optimal shrinkage intensities given in (\ref{ea}) and (\ref{eb}) with the corresponding ones calculated for the population covariance matrix given in Bodnar et al. (2013) then we conclude that the corresponding estimators are different although the structure remains somewhat similar. In the simulation study of Section 5 both estimators for the precision matrix are compared with each other and it is shown that it is better to shrink the inverse sample covariance matrix itself than to shrink the sample covariance matrix and then to invert it.

\subsection{Choice of $\boldsymbol{\Pi}_0$.}
The next question is the choice of the nonrandom target matrix $\boldsymbol{\Pi}_0$ which should be positive definite with uniformly bounded trace norm. Unfortunately, the answer on this question is not yet found because the choice of the target matrix is equivalent to the choice of the hyperparameter for the prior distribution of $\bSigma^{-1}_n$. This problem is well-known in Bayesian statistics. The application of different priors leads to different results. So it is very important to choose the one which works fine in most cases. The naive one is $\boldsymbol{\Pi}_0=\frac{1}{p}\bI$ where $\bI$ is the identity matrix. Obviously, the oracle OLSE estimator has the prior matrix as the true precision matrix $\bSigma^{-1}_n$ and is a consistent estimator for the precision matrix as shown in Proposition 4.1. Moreover, in the next section we show how the prior information on the spectrum of the precision matrix can significantly improve the OLSE estimator.

Consider the OLSE estimator as a matrix function $\widehat{\boldsymbol{\Pi}}_{OLSE}(\boldsymbol{\Pi}_0):M_p\rightarrow \tilde{M}_p$, where $M_p$ is the space of $p$-dimensional positive definite symmetric matrices and $\tilde{M}_p$ is the corresponding space of the $p$-dimensional positive definite symmetric random matrices. In the following proposition we present some properties of the OLSE estimator as a function of the prior matrix $\boldsymbol{\Pi}_0$.
\begin{proposition}
For the OLSE estimator $\widehat{\boldsymbol{\Pi}}_{OLSE}(\boldsymbol{\Pi}_0)$ it holds that
\begin{enumerate}[label=\roman{*}).]
\item the function $\widehat{\boldsymbol{\Pi}}_{OLSE}(\boldsymbol{\Pi}_0)$ is scale invariant, i.e. for arbitrary $\sigma>0$ $\widehat{\boldsymbol{\Pi}}_{OLSE}(\sigma\boldsymbol{\Pi}_0)=\widehat{\boldsymbol{\Pi}}_{OLSE}(\boldsymbol{\Pi}_0)$.
\vspace{3mm}
\item $\widehat{\boldsymbol{\Pi}}_{OLSE}(1/p\bI)$ is a consistent estimator for the precision matrix $\bSigma^{-1}_n=\sigma\bI$ for arbitrary $\sigma>0$ and $c\in(0, +\infty)$.
\vspace{3mm}
\item $\widehat{\boldsymbol{\Pi}}_{OLSE}(\bSigma_n^{-1})$ is a consistent estimator for the precision matrix $\bSigma^{-1}_n$.
\end{enumerate}
\end{proposition}

\section{Simulation Study} In this section we provide a Monte Carlo simulation study to investigate the performance of the suggested OLSE estimator for the precision matrix.

Before we proceed, two benchmark estimators for the precision matrix used in the simulations are introduced. The first one was already mentioned in Section 4 and it is called the OLSE estimator for the population covariance matrix. This estimator has been suggested by Bodnar et al. (2013) and is a generalization of the well-conditioned estimator proposed by Ledoit and Wolf (2004). The second estimator is the so-called nonlinear shrinkage estimator developed by Ledoit and Wolf (2012) which seems to be the best one up to now.

The nonlinear shrinkage estimator is based on the class of rotation-equivariant estimators for the precision matrix, i.e., on the estimators which have the same eigenvectors as the sample covariance matrix. Ledoit and Wolf (2012) proposed the following oracle equivariant (EV) estimator for $\bSigma^{-1}_n$
\begin{equation}\label{est}
\widehat{\boldsymbol{\Pi}}_{EV}=\mathbf{U}\mathbf{A}^*\mathbf{U}^{\prime}~~\text{with}~~\mathbf{A}^*=\text{diag}\{\mathbf{U}^{\prime}\bSigma_n^{-1}\mathbf{U}\}=\{a^*_i\}^{p}_{i=1}\,,
\end{equation}
where the matrix $\mathbf{U}=(\mathbf{u}_1,\ldots,\mathbf{u}_p)$ is the eigenmatrix (matrix of the orthonormal eigenvectors) of the sample covariance matrix $\bS_n$ and the diagonal matrix $\mathbf{A}^*$ is the unique minimizer of the Frobenius loss, i.e.,
\begin{align}
\mathbf{A}^*&=\text{argmin}\left(||\mathbf{U}\mathbf{A}\mathbf{U}^{\prime}-\bSigma^{-1}_n||^2_F\right)~~\text{for a diagonal matrix}~\mathbf{A}\label{A}\,.
\end{align}

The second benchmark in our study is the inverse of the OLSE estimator for the population covariance matrix developed by Bodnar et al. (2013). Recall, that the optimal linear shrinkage estimator (OLSE) for the covariance matrix $\bSigma_n$ for $c\in(0, +\infty)$ is given by
 \begin{equation}\label{olsep}
\widehat{\bSigma}_{OLSE}=\tilde{\alpha}^*\bS_n+\tilde{\beta}^*\bSigma_0~~\text{with}~||\bSigma_0||_{tr}\leq M\,,
\end{equation}
where $\bSigma_0$ is the positive definite symmetric target matrix,
\begin{equation}\label{a1}
 \tilde{\alpha}^*=1-\dfrac{\dfrac{1}{n}||\bS_n||^2_{tr}||\bSigma_0||^2_F}{||\bS_n||^2_F||\bSigma_0||^2_F-\bigl(\text{tr}(\bS_n\bSigma_0)\bigr)^2}~~\text{and}~~\tilde{\beta}^*=\dfrac{\text{tr}(\bS_n\bSigma_0)}{||\bSigma_0||^2_F}\left(1-\tilde{\alpha}^*\right)\,.
 \end{equation}
Bodnar et al. (2013) proved that their $\widehat{\bSigma}_{OLSE}$ estimator possesses asymptotically almost surely the smallest Frobenius loss over all linear shrinkage estimators. Moreover, they showed that if the target matrix is $\bSigma_0=1/p\bI$ then the estimator $\widehat{\bSigma}_{OLSE}$ is asymptotically equivalent to the linear shrinkage estimator proposed by Ledoit and Wolf (2004).
Of course, in order to compare this estimator with the suggested OLSE estimator for the precision matrix $\widehat{\boldsymbol{\Pi}}_{OLSE}$ from (\ref{olse}) we have to invert $\widehat{\bSigma}_{OLSE}$.


Next, we compare the performance of the estimators $\widehat{\boldsymbol{\Pi}}_{OLSE}$, $\widehat{\bSigma}^{-1}_{OLSE}$ and $\widehat{\boldsymbol{\Pi}}_{EV}$ given in (\ref{olse}), (\ref{olsep}) and (\ref{est}), respectively. As a performance measure we take the PRIAL (Percentage Relative Improvement in Average Loss) presented in Ledoit and Wolf (2012). For an arbitrary estimator of the precision matrix, $\widehat{\mathbf{M}}$, the PRIAL is defined by
\begin{equation}\label{PRIAL}
 \text{PRIAL}(\widehat{\mathbf{M}})=\left(1-\dfrac{E||\widehat{\mathbf{M}}-\bSigma^{-1}_n||^2_F}{E||\bS^{-1}_n-\bSigma^{-1}_n||_F^2}\right)\cdot100\%\,.
 \end{equation}
Thus, by definition (\ref{PRIAL}), PRIAL($\bS^{-1}_n$) is equal to zero and PRIAL($\bSigma^{-1}_n$) is equal to 100$\%$.

In our simulations, without loss of generality, we take $\bSigma_n$ as a diagonal matrix and separate its spectrum in three parts with $20\%$ of the eigenvalues equal to $1$, $40\%$ equal to $3$ and $40\%$ equal to $10$. In terms of the corresponding cumulative distribution function of the eigenvalues of $\bSigma_n$ (cf. Section 2) it holds that
 \begin{equation}\label{sigmaH}
 H^{\bSigma_n}_n(t)=1/5\delta_{[1,\hspace{1mm}\infty)}(t)+2/5\delta_{[3,\hspace{1mm}\infty)}(t)+2/5\delta_{[10,\hspace{1mm}\infty)}(t)\,,
 \end{equation}
 where $\delta$ is the Dirac delta function. Doing so we leave the structure of population covariance matrix unchanged for all dimensions $p$. Note that the same structure of the population covariance matrix was also used in the simulation study by Ledoit and Wolf (2012).

The next step is to choose the shrinkage target matrices $\bSigma_0$ and $\boldsymbol{\Pi}_0$ for the OLSE estimators $\widehat{\bSigma}^{-1}_{OLSE}$ and $\widehat{\boldsymbol{\Pi}}_{OLSE}$, respectively. As the first prior target we choose a simplest one $\boldsymbol{\Pi}_0=\bSigma_0=1/p\bI$. In the choice of the second prior we want to concentrate on the prior information about the spectrum of the covariance (precision) matrix. Thus, we assume that the information about the spectrum separation of the population covariance matrix is available which is separated in three blocks in the relation $1/5:2/5:2/5$ (see, equality (\ref{sigmaH})). No other information is taken into account. The diagonal elements of the prior matrix $\bSigma_0$ are chosen, without loss of generality, to be $1$, $2$ and $4$. In terms of the cumulative distribution function of $\bSigma_0$ it holds
\begin{equation}\label{sigmaH0}
 H^{\bSigma_0}(t)=1/5\delta_{[1,\hspace{1mm}\infty)}(t)+2/5\delta_{[2,\hspace{1mm}\infty)}(t)+2/5\delta_{[4,\hspace{1mm}\infty)}(t)\,.
 \end{equation}
Note that the eigenvalues $1$, $2$ and $4$ are far away from the real population eigenvalues presented in (\ref{sigmaH}) and are chosen with the aim to depart from the naive prior $1/p\bI$. Indeed, further we show that for the choice of the prior $\bSigma_0$ the knowledge of the spectrum separation is more important than the values on the diagonal itself. It is noted that the corresponding target matrix $\boldsymbol{\Pi}_0$ is chosen as the inverse of $\bSigma_0$, i.e.
\begin{equation}\label{sigmaHPi0}
 H^{\boldsymbol{\Pi}_0}(t)=1/5\delta_{[1,\hspace{1mm}\infty)}(t)+2/5\delta_{[1/2,\hspace{1mm}\infty)}(t)+2/5\delta_{[1/4,\hspace{1mm}\infty)}(t)\,.
 \end{equation}

In Figure 1 we present the first simulation results for the normally distributed data when $c=1/3$. The oracle estimators are presented as solid lines while the corresponding bona fide estimators are dashed lines. For all bona fide estimators we observe a fast convergence rate to their corresponding oracles, i.e. for the dimension $p\geq50$ they all converge in PRIAL.
It is remarkable that the bona fide OLSE estimator for the precision matrix $\widehat{\boldsymbol{\Pi}}_{OLSE}$ with the naive prior $\boldsymbol{\Pi}_0=1/p\bI$ dominates the corresponding inverted OLSE estimator $\widehat{\bSigma}^{-1}_{OLSE}$ for all chosen priors. We know that the OLSE estimator $\widehat{\bSigma}_{OLSE}$ with prior $1/p\bI$ is asymptotic equivalent to the Ledoit-Wolf linear shrinkage (see, Bodnar et al. (2013)). The most stunning result we observe is for the bona fide OLSE estimator $\widehat{\boldsymbol{\Pi}}_{OLSE}$ with the prior $\boldsymbol{\Pi}_0=\bSigma^{-1}_0$ which contains the information on the spectrum separation. It dominates the oracle equivariant estimator $\widehat{\boldsymbol{\Pi}}_{EV}$ and thus it dominates the nonlinear shrinkage estimator $\widehat{\boldsymbol{\Pi}}_{LW}$. Next, we compare the results of Figure 1 with those of Figure 6 presented in Ledoit and Wolf (2012). It appears that the inverse of the nonlinear shrinkage estimator for the covariance matrix is on the similar level as the OLSE estimator $\widehat{\boldsymbol{\Pi}}_{OLSE}$ with the naive prior $\boldsymbol{\Pi}_0=1/p\bI$. This result ensures that for the estimation of the precision matrix the simple OLSE estimator with the naive prior $1/p\bI$ is a great alternative to the inverted nonlinear shrinkage estimator for the population covariance matrix.
\begin{figure}[H]
\centerline{\includegraphics[scale=0.45]{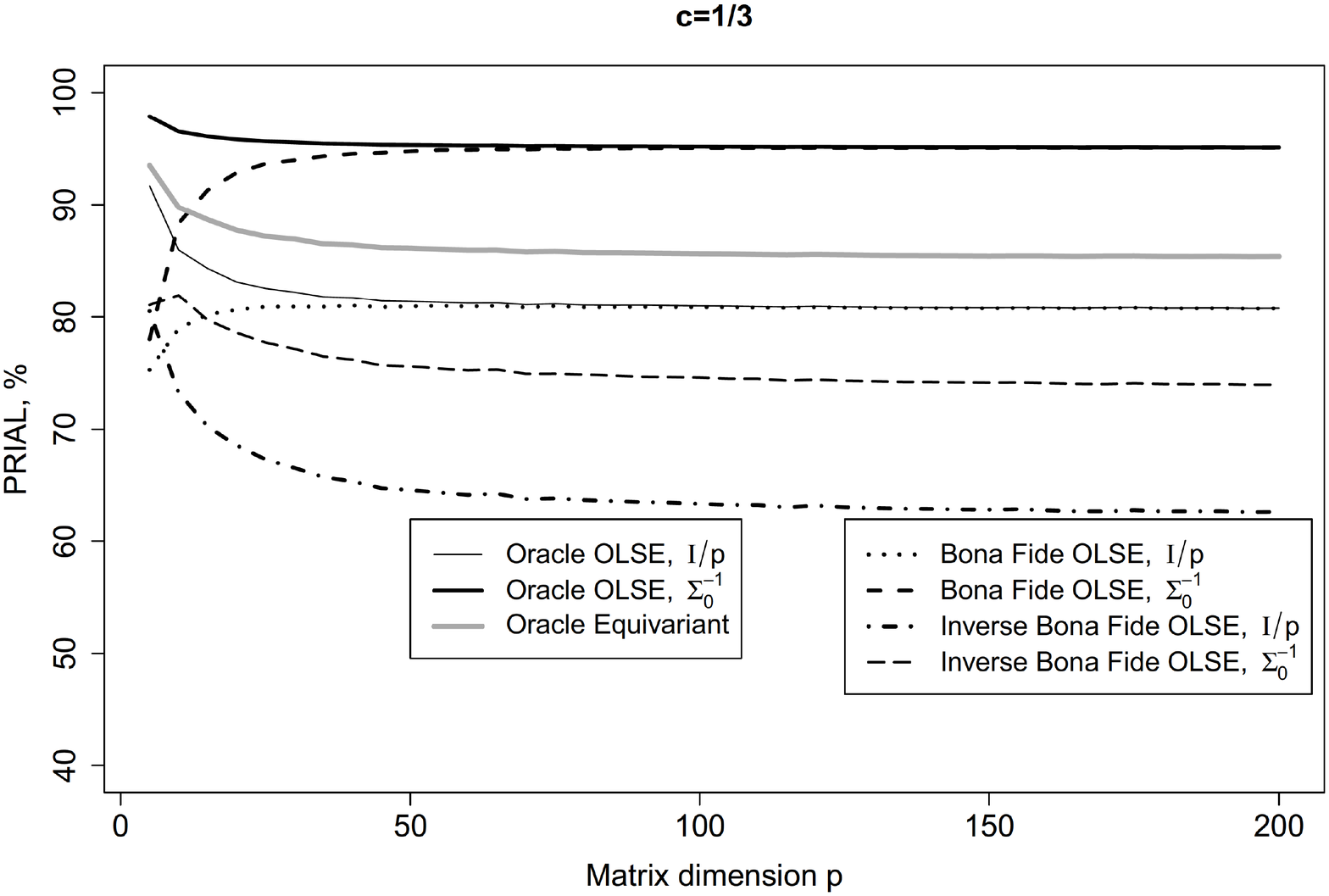}}
\caption{PRIALs for the oracle and bona fide estimator $\widehat{\boldsymbol{\Pi}}_{OLSE}$ with the prior $\boldsymbol{\Pi}_0=\bSigma^{-1}_0$ and with the $\boldsymbol{\Pi}_0=1/p\bI$, the inverse of the bona fide estimator $\widehat{\bSigma}_{OLSE}$ with the prior $\bSigma_0$ and wtih the prior $1/p\bI$, and the oracle equivariant estimator $\widehat{\boldsymbol{\Pi}}_{EV}$ given in (\ref{est}) for $p=5k,k\in\{1,\ldots,40\}$, $c=1/3$. The results are based on $1000$ independent realizations.}
\label{Fig:1}
\end{figure}
More interesting results are present in Figure 2. Here we show that our argument about different priors holds true: the information about the spectrum separation is sufficient to construct a dominating OLSE estimator for the precision matrix. We take 5 different prior matrices $\bSigma_0$, which contain the information about the spectrum separation of $\bSigma^{-1}_n$ (\ref{sigmaH}). In terms of the corresponding (e.d.f.) it holds
\begin{align*}
H^{\bSigma^{(1)}_0}(t)&=1/5\delta_{[1,\hspace{1mm}\infty)}(t)+2/5\delta_{[5,\hspace{1mm}\infty)}(t)+2/5\delta_{[10,\hspace{1mm}\infty)}(t)\\
H^{\bSigma^{(2)}_0}(t)&=1/5\delta_{[1,\hspace{1mm}\infty)}(t)+2/5\delta_{[2,\hspace{1mm}\infty)}(t)+2/5\delta_{[4,\hspace{1mm}\infty)}(t)\\
H^{\bSigma^{(3)}_0}(t)&=1/5\delta_{[1,\hspace{1mm}\infty)}(t)+2/5\delta_{[2,\hspace{1mm}\infty)}(t)+2/5\delta_{[60,\hspace{1mm}\infty)}(t)\\
H^{\bSigma^{(4)}_0}(t)&=1/5\delta_{[0.1,\hspace{1mm}\infty)}(t)+2/5\delta_{[1,\hspace{1mm}\infty)}(t)+2/5\delta_{[1000,\hspace{1mm}\infty)}(t)\\
H^{\bSigma^{(5)}_0}(t)&=1/5\delta_{[0.1,\hspace{1mm}\infty)}(t)+2/5\delta_{[0.5,\hspace{1mm}\infty)}(t)+2/5\delta_{[1,\hspace{1mm}\infty)}(t)\,.
\end{align*}

\begin{figure}[H]
\centerline{\includegraphics[scale=0.45]{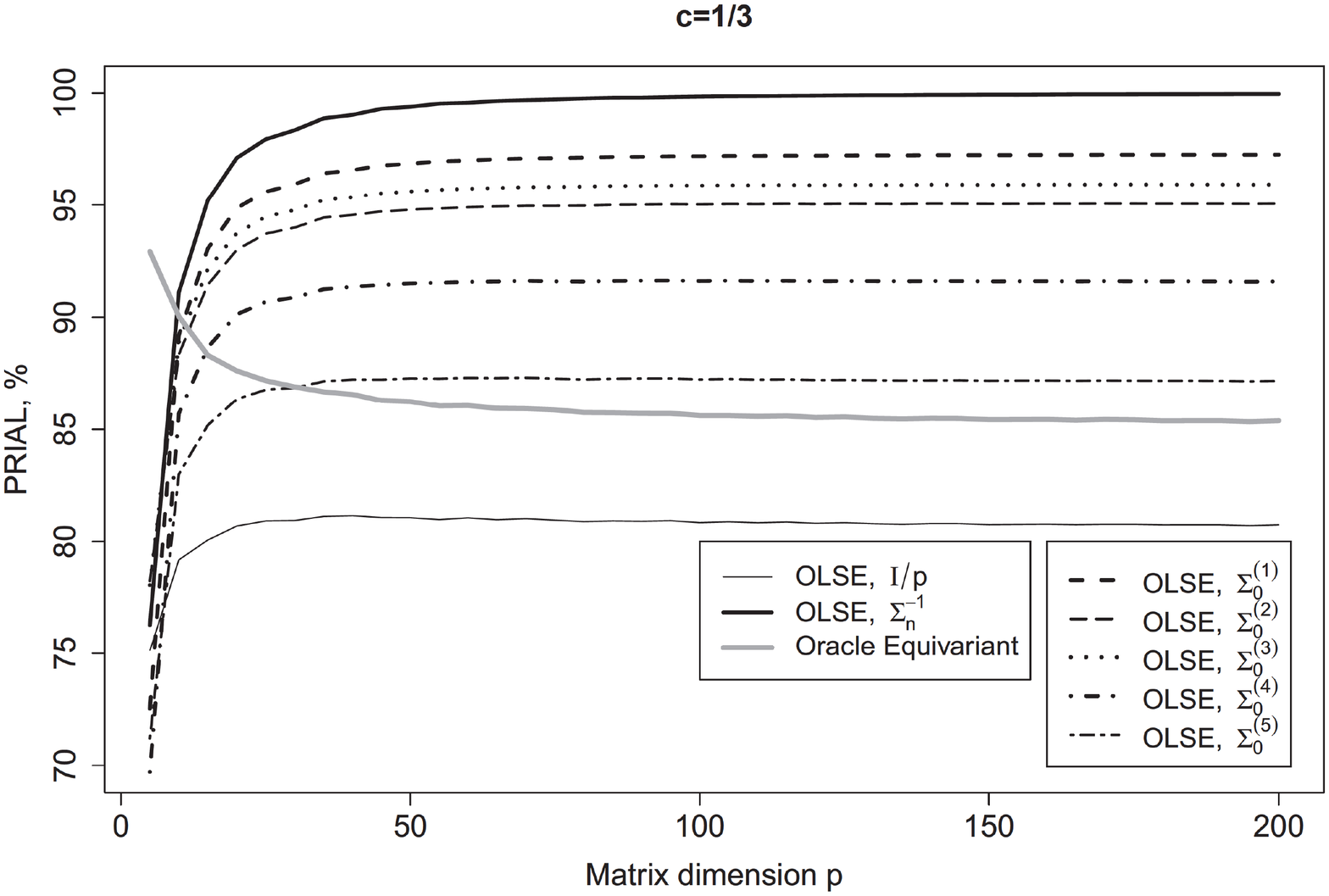}}
\caption{PRIALs for the bona fide estimator $\widehat{\boldsymbol{\Pi}}_{OLSE}$ with different priors which contain information about spectrum separation, with the naive prior $1/p\bI$, with the prior which equals the true precision matrix $\bSigma^{-1}_n$, and the oracle equivariant estimator $\widehat{\boldsymbol{\Pi}}_{EV}$ given in (\ref{est}) for $p=5k,k\in\{1,\ldots,40\}$, $c=1/3$. The results are based on $1000$ independent realizations.}
\label{Fig:2}
\end{figure}

The corresponding $\boldsymbol{\Pi}^{(i)}_0$ is constructed by inverting the prior $\bSigma^{(i)}_0$ given above.
Figure 2 shows that all of the five OLSE estimators for the precision matrix with differently chosen priors dominate the oracle equivariant estimator $\widehat{\boldsymbol{\Pi}}_{EV}$.
In Figure 2 we also present the bona fide OLSE estimator with the prior $\bSigma^{-1}_n$ which obviously converges to the precision matrix in PRIAL and forms the upper bound of all estimators (see, Proposition 4.1); the lower bound builds the naive prior $1/p\bI$. As expected, the best prior is $\bSigma^{(1)}_0$ which contains two eigenvalues from the population covariance matrix. The worst one is $\bSigma^{(5)}_0$, which is also expected, because the distance between the largest and the smallest eigenvalues is small and thus the blocks of the spectrum are near to each other. The results obtained do not mean that there is no prior which gives worse results than the equivariant estimator. However, they show that if we know the spectrum separation of $\bSigma^{-1}_n$ or at least can estimate it consistently then, without much effort, we are able to construct a dominating OLSE estimator for the precision matrix. Thus, the suitably chosen prior can significantly improve the estimator of the precision matrix given in (\ref{olse}). Nevertheless, it would be a challenging task to prove analytically the findings obtained via the simulations. Note that the similar results arise when we assume that the population covariance matrix is not a diagonal one.

In Figure 3 we present the results of simulations under the normal distribution for $c=1/2$ and $c=0.8$. For $c=1/2$ we observe a better overall performance for all of the considered estimators for the precision matrix. The OLSE estimator (\ref{olse}) with the prior given in (\ref{sigmaHPi0}) is again the best one and it is almost consistent. On the second place we put the oracle equivariant estimator $\widehat{\boldsymbol{\Pi}}_{EV}$ while on the third place both $\widehat{\boldsymbol{\Pi}}_{OLSE}(1/p\bI)$ and $\widehat{\bSigma}^{-1}_{OLSE}$ from (\ref{olsep}) with the prior given in (\ref{sigmaH0}) are placed. The worst one in the case $c=1/2$ is the inverse of the OLSE estimator  $\widehat{\bSigma}_{OLSE}(1/p\bI)$.

Figure 3 for $c=0.8$ shows that all estimators provide the superior performance but the best one is again $\widehat{\boldsymbol{\Pi}}_{OLSE}$ with prior (\ref{sigmaHPi0}) which is consistent in this case. The second place belongs to $\widehat{\bSigma}^{-1}_{OLSE}$ with prior (\ref{sigmaH0}) and the equivariant estimator $\widehat{\boldsymbol{\Pi}}_{EV}$ from (\ref{est}) which surprisingly converge to each other. The last ones are the OLSE estimators  $\widehat{\boldsymbol{\Pi}}_{OLSE}(1/p\bI)$ and $\widehat{\bSigma}_{OLSE}(1/p\bI)$, respectively. From the last simulation for $c=0.8$ it can be observed that all of the considered estimators are near to the true precision matrix.

Figure 4 is dedicated to the case when the underlying distribution departs from the normal one. Here we consider the $t$-distribution with 10 degrees of freedom. The estimators and their priors are the same as used in Figure 1. For $c=1/3$ we observe that the overall performance of all considered estimators is even better as in the case of the normal distribution in Figure 1. As usual, the OLSE estimator with the prior (\ref{sigmaHPi0}) is ranked first. It is remarkable that the suggested bona fide OLSE estimators converge much slower to their oracles. In contrast to Figure 1 where the convergence was very fast ($p\geq50$), Figure 4 ensures the convergence for $p\geq500$. It seems the convergence rate is influenced by heavy tails, this shows that the heavier the tails are the slower is the convergence of the bona fide OLSE estimators to their oracles. The second place belongs to the oracle equivariant estimator (EV) which seems to be quite robust to the presence of heavy tail. Note that in the third place are the inverse bona fide estimator  $\widehat{\bSigma}^{-1}_{OLSE}$ with the prior (\ref{sigmaH0}) and the bona fide $\widehat{\boldsymbol{\Pi}}_{OLSE}(1/p\bI)$ with its oracle. It seems that the inverse bona fide estimator $\widehat{\bSigma}^{-1}_{OLSE}$ and the bona fide estimator $\widehat{\boldsymbol{\Pi}}_{OLSE}(1/p\bI)$ converge to the same oracle. The last one, as usual, is the inverse bona fide estimator $\widehat{\boldsymbol{\Sigma}}^{-1}_{OLSE}(1/p\bI)$ which is asymptotically equivalent to the inverse linear shrinkage estimator proposed by Ledoit and Wolf (2004). As a result, even in the non-normal case the proposed OLSE estimator for the precision matrix shows remarkable stability and robustness. The simulation results for $c=1/2$ and $c=0.8$ are very similar to those obtained for the normal distribution in Figure 3. Only the convergence of the bona fide estimators to their oracles is slower.

\begin{figure}[H]
\begin{tabular}{c}
\centerline{\includegraphics[scale=0.423]{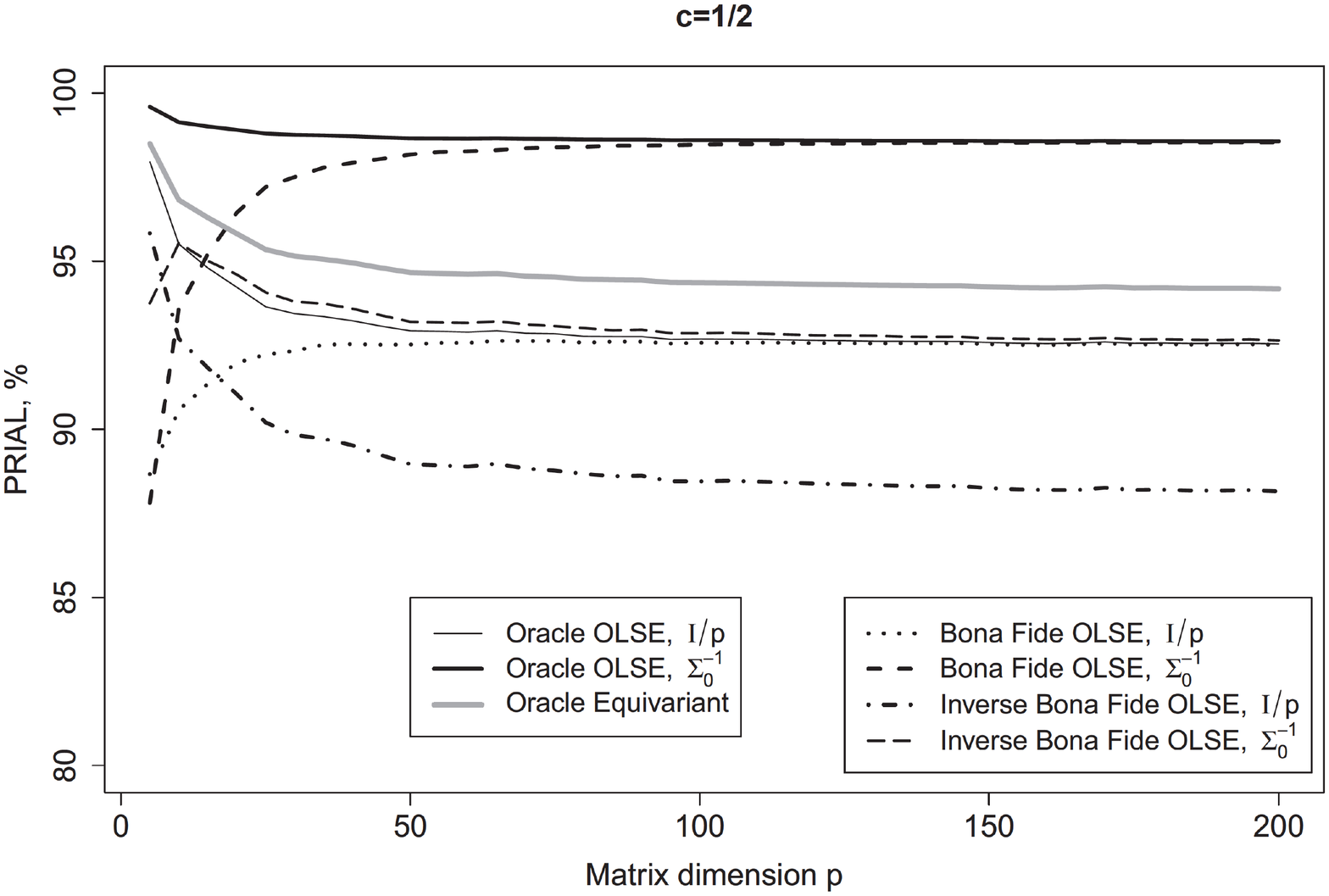}}
\\
\centerline{\includegraphics[scale=0.423]{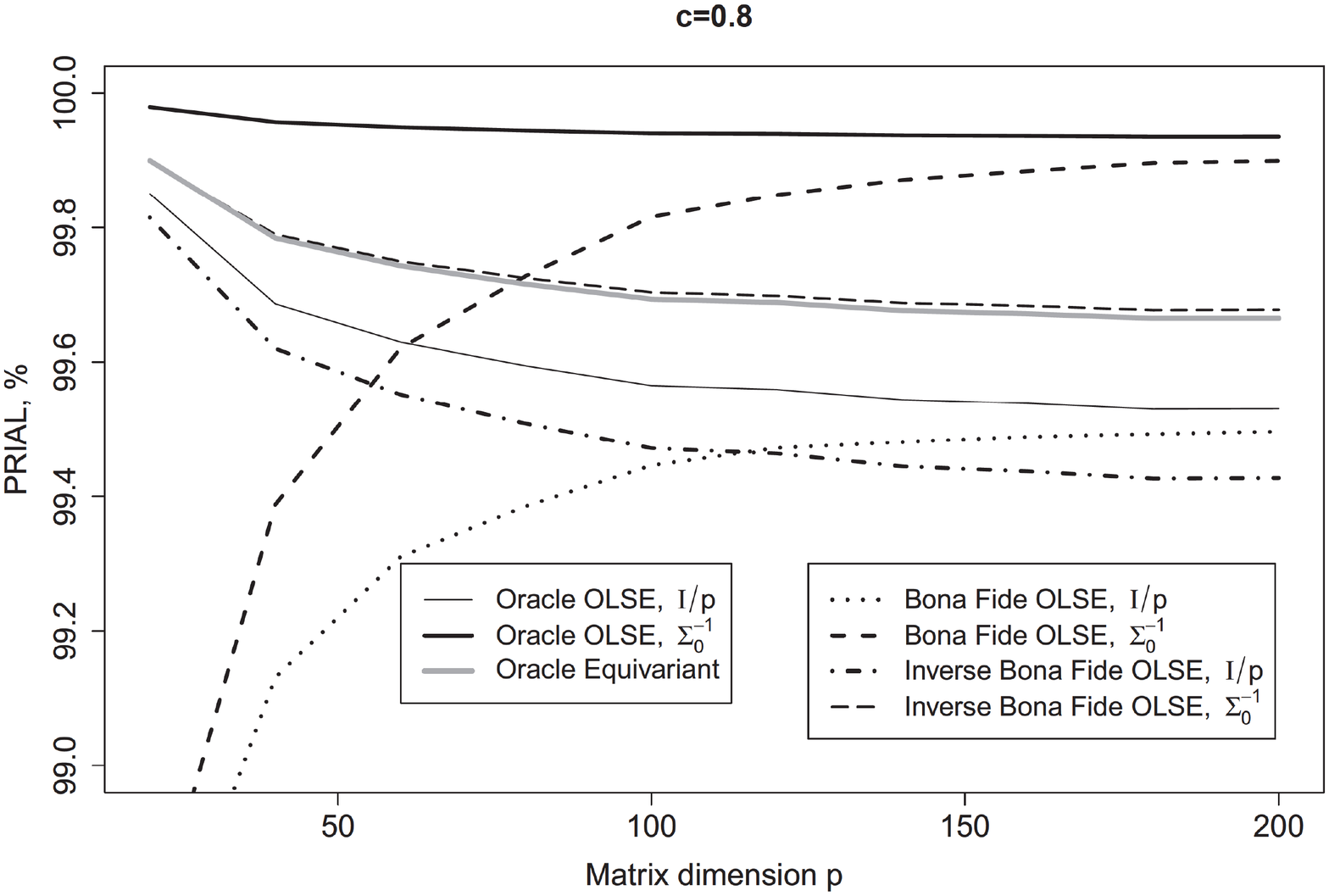}}
\end{tabular}
\caption{PRIALs for the oracle and the bona fide estimator $\widehat{\boldsymbol{\Pi}}_{OLSE}$ with the prior $\boldsymbol{\Pi}_0=\bSigma^{-1}_0$ and with the prior $\boldsymbol{\Pi}_0=1/p\bI$, the inverse of the bona fide estimator $\widehat{\bSigma}_{OLSE}$ with the prior $\bSigma_0$ and with  the prior $1/p\bI$, and the oracle equivariant estimator $\widehat{\boldsymbol{\Pi}}_{EV}$ given in (\ref{est}). We put $p=5k,k\in\{1,\ldots,40\}$ for $c=1/2$ (upper figure) and $p=20k,k\in\{1,\ldots,10\}$ for $c=0.8$ (lower figure). The results are based on $1000$ independent realizations.}
\label{Fig:3}
\end{figure}

\begin{figure}[H]
\vspace{-0.4cm}
\centerline{\includegraphics[scale=0.423]{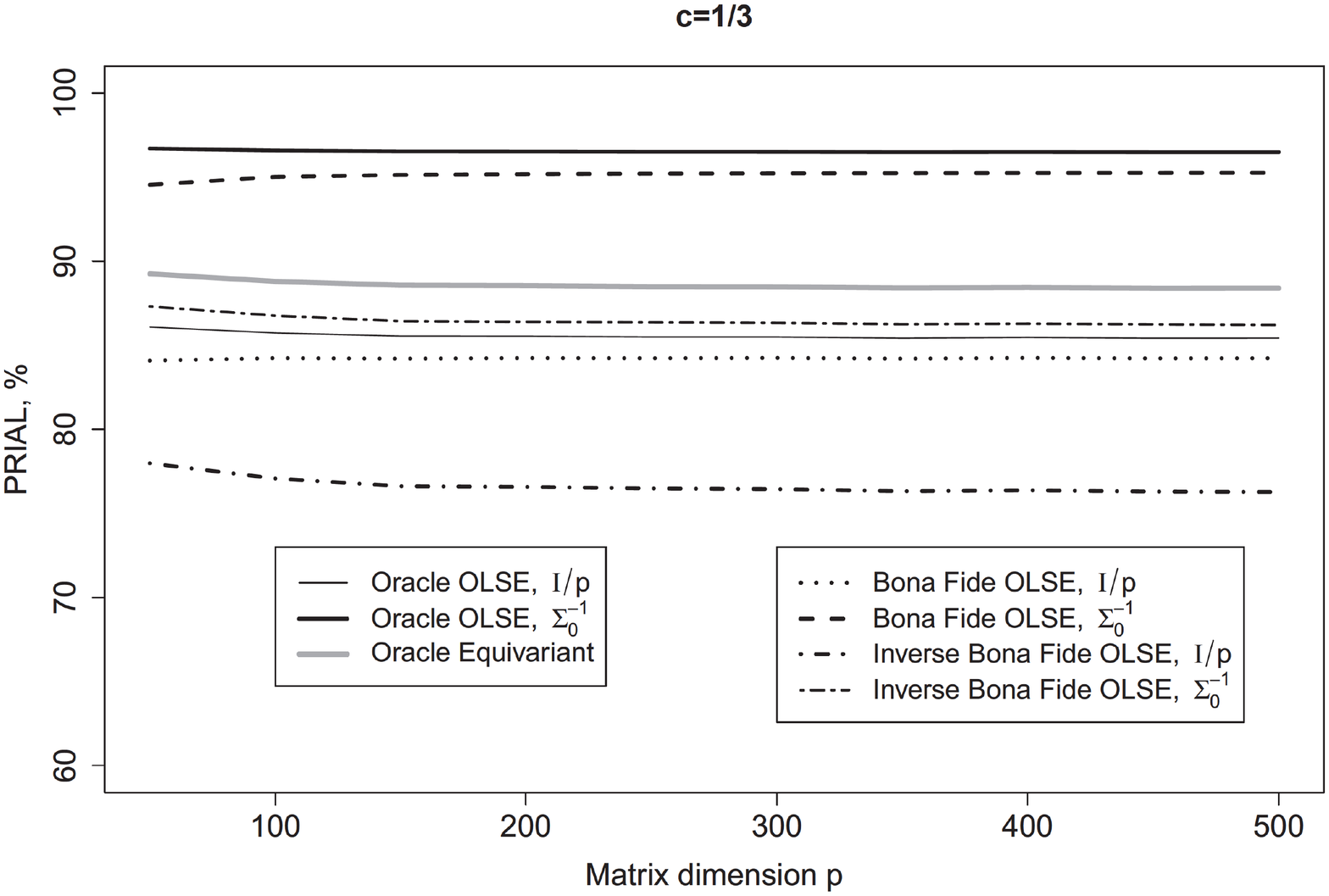}}
\caption{PRIALs for the oracle and the bona fide estimator $\widehat{\boldsymbol{\Pi}}_{OLSE}$ with the prior $\boldsymbol{\Pi}_0=\bSigma^{-1}_0$ and with the prior $\boldsymbol{\Pi}_0=1/p\bI$, the inverse of the bona fide $\widehat{\bSigma}_{OLSE}$ with the prior $\bSigma_0$ and with the prior $1/p\bI$, and the oracle equivariant estimator $\widehat{\boldsymbol{\Pi}}_{EV}$ given in (\ref{est}) for $p=50k,k\in\{1,\ldots,10\}$, $c=1/3$. The data are generated from the $t$-distribution with $10$ degrees of freedom. The results are based on $1000$ independent realizations.}
\label{Fig:4}
\end{figure}

\begin{figure}[H]
\vspace{-0.4cm}
\centerline{\includegraphics[scale=0.423]{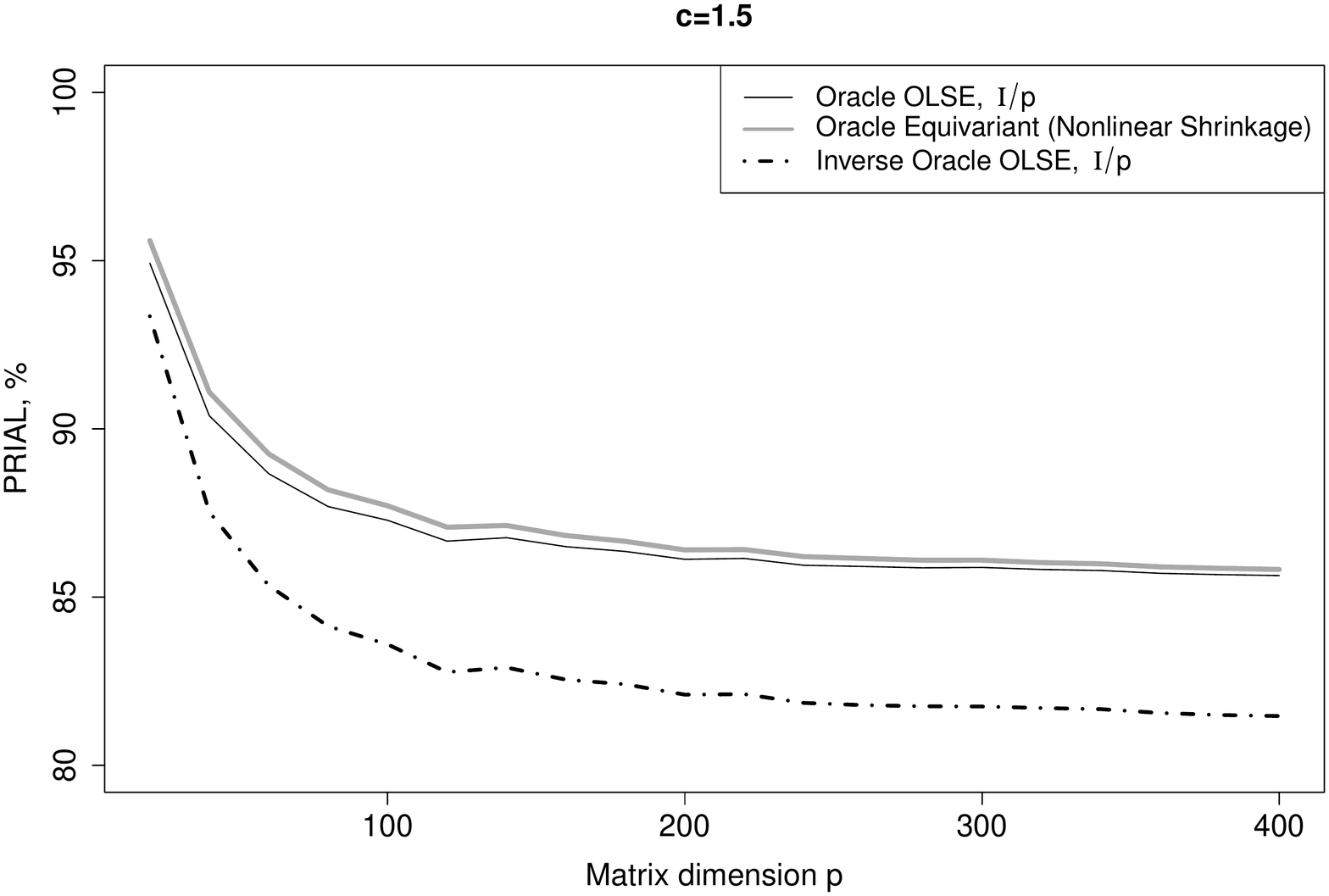}}
\caption{PRIALs for the oracle estimator $\widehat{\boldsymbol{\Pi}}_{OLSE}$ with the prior $\boldsymbol{\Pi}_0=1/p\bI$, the inverse of the bona fide $\widehat{\bSigma}_{OLSE}$ with the prior $1/p\bI$, and the oracle equivariant estimator $\widehat{\boldsymbol{\Pi}}_{EV}$ given in (\ref{est}) for $p=20k,k\in\{1,\ldots,20\}$, $c=1.5$. The data are generated from the normal distribution. The results are based on $1000$ independent realizations.}
\label{Fig:5}
\end{figure}

At last, Figure 5 contains the example when $c=1.5$ for the normally distributed data. Here we use the oracle OLSE estimator proposed in Section 3.2 and compare it with the oracle nonlinear shrinkage estimator, the oracle inverse OLSE and the generalized inverse of the sample covariance matrix $\bS^+_n$. The PRIAL is defined in the similar way as (\ref{PRIAL}), only the matrix $\bS^{-1}$ is changed to the generalized inverse $\bS^+_n$. Similarly, in case of $\bS^+_n$ the PRIAL is equal to zero. Surprisingly, the nonlinear shrinkage estimator and the suggested oracle OLSE estimator are converging to each other. They show the performance at roughly $87\%$ level while the bona fide inverse OLSE estimator at ca. $83\%$. The same convergence is detected for $c=2$\footnote{This result is not presented in the paper and it is available from authors on request.} but the overall performance of all estimators has decreased. This means that if we could consistently estimate the shrinkage intensities given by (\ref{talfa1}) and (\ref{tbeta1}), then the resulting bona fide OLSE estimator would probably converge to the nonlinear shrinkage proposed by Ledoit and Wolf (2012). This task is not an easy one and leaves the place for the future research.

As a result, the simulation results as well as the theoretical findings show that the OLSE estimator $\widehat{\boldsymbol{\Pi}}_{OLSE}(\boldsymbol{\Pi})$ is a great alternative not only to the sample estimator and to the inverted linear shrinkage estimator proposed by Ledoit and Wolf (2004) and generalized by Bodnar et. al (2013) but also to the nonlinear shrinkage estimator recently suggested by Ledoit and Wolf (2012). The case of $c>1$ is even more important for the practical purposes but it seems to be more difficult to handle analytically. This can be done in an efficient way if the population covariance matrix is a multiple of identity. In general case a good alternative would be the inverse OLSE estimator given in (\ref{olsep}), but it is not optimal for the precision matrix. This point will be treated in future research.

\section{Summary} In this paper the problem of the estimation of the precision matrix for large dimensional data is considered. Our particular interest is the case when both the dimension of the precision matrix $p\rightarrow\infty$ and the sample size $n\rightarrow\infty$ so that $p/n\rightarrow c\in (0, +\infty)$. Using the results from the random matrix theory and the linear shrinkage technique we develop an estimator for the precision matrix which is distribution-free, robust and possesses almost surely the smallest Frobenius loss asymptotically. In particular, we prove that the Frobenius norms of the inverse and of the generalized inverse sample covariance matrices as well as of the optimal shrinkage intensities tend to the nonrandom quantities under high dimensional asymptotics. In order to get the optimal linear shrinkage estimator for the precision matrix we estimate the unknown quantities consistently. The performance of the suggested OLSE estimator is compared with other known estimators for the precision matrix via the simulation study. The OLSE estimator shows significant improvements in the presence of the prior information about the spectrum separation of the population precision matrix.

\section*{Acknowledgments}
The authors would like to thank Dr. C. Guillaume for his valuable comments on the choice of target matrix and the practical insights of the shrinkage estimation.

\section{Appendix} Here the proofs of the theorems are given.\\

\noindent\textbf{Proof of Theorem 3.1}.
The proof of the theorem is based on the Marchenko-Pastur theorem proved by Silverstein (1995).
\begin{theorem}{\textbf{[Silverstein (1995)]}} Assume that on the common probability space assumption (A1) is satisfied for $\frac{p}{n}\rightarrow c\in(0,+\infty)$ as $n\rightarrow\infty$. Then almost surely $F_n(t)\stackrel{a.s.}{\Rightarrow}F(t)$ as $n\rightarrow\infty$. Moreover, the Stieltjes transform of $F$ satisfies the following equation
\begin{equation}\label{MP}
m_F(z)=\int\limits_{-\infty}^{+\infty}\dfrac{1}{\tau(1-c-czm_F(z))-z}dH(\tau)\,,
\end{equation}
in the sense that $m_F(z)$ is the unique solution of (\ref{MP}) for all $z\in\mathbbm{C}^+$.
\end{theorem}

Consider the asymptotics of the quantity
\begin{equation}\label{Frobs}
\dfrac{1}{p}\text{tr}\left(\bS_n^{-2}\right)=\left.\dfrac{\partial}{\partial z}\dfrac{1}{p}\text{tr}\left[(\bS_n-z\bI)^{-1}\right]\right|_{z=0}=\left.\dfrac{\partial}{\partial z}m_{F_n}(z)\right|_{z=0}\,.
\end{equation}
Using Theorem 7.1 we note that $m_{F_n}(z)$ tends almost surely to a nonrandom limit function $m _{F}(z)$ which is the unique solution of the MP equation (\ref{MP})
First, we show that the limit $z\rightarrow0^+$ can be taken under the integral sign in (\ref{MP}). Let $\underline{m}(z)=-\dfrac{1-c}{z}+cm_F(z)$ then using the assumption (A3) we rewrite equality (\ref{MP}) in the following way
\begin{align}\label{rew1}
m_F(z)&=\int\limits_{-\infty}^{+\infty}\dfrac{-1}{z(\tau\underline{m}(z)+1)}dH(\tau)=\int\limits_{h_0}^{h_1}\dfrac{-1}{z(\tau\underline{m}(z)+1)}dH(\tau)
=\int\limits_{h_0}^{h_1}\dfrac{s(z,\tau)}{\tau}dH(\tau)\,,
\end{align}
where the function $s(z,\tau)=\dfrac{-\tau}{z(\tau\underline{m}(z)+1)}$ is the Stieltjes transform of a positive measure on $\mathbbm{R}^+$ with total mass $\tau$ (see, Paul and Silverstein (2009)).
Using the properties of Stieltjes transform $\overline{s(z,\tau)}=s(\overline{z},\tau)$ and $\lim\limits_{z\rightarrow0^+}s(z,\tau)=\lim\limits_{z\rightarrow0^+}s(\overline{z},\tau)$ we get
\begin{equation}\label{lmods}
\gamma^2(c)\equiv\lim\limits_{z\rightarrow0^+}|s(z,\tau)|^2=\lim\limits_{z\rightarrow0^+}s(z,\tau)s(\overline{z},\tau)=\left(\lim\limits_{z\rightarrow0^+}s(z,\tau)\right)^2=\left(\dfrac{1}{1-c}\right)^2\,.
\end{equation}
The last equality in (\ref{lmods}) follows from boundedness of $m_F(0)=\lim\limits_{z\rightarrow0^+}m_F(z)$ for $c<1$ (see, Silverstein and Choi (1995)).

Let $\mathcal{B}$ be a compact ball around $0$. Then, the application of the inequality (cf. Silverstein (2009))
\begin{equation}
0\leq\left|s(z,\tau)\right|\leq\dfrac{1}{\mathbf{Im}(z)}\,.
\end{equation}
leads to


\begin{equation}\label{modf}
0\leq\left|\dfrac{s(z,\tau)}{\tau}\right|\leq\dfrac{\gamma(c)}{h_0}~~\forall z\in\mathcal{B}\,.
\end{equation}
Hence, using the inequality (\ref{modf}) together with (\ref{rew1}) and the dominated convergence theorem we conclude that the limit $z\rightarrow0^+$ can be moved under the integral sign in (\ref{MP}). This fact together with $m_F(0)<\infty$ implies that
\begin{equation}\label{lim0}
m_F(0)=\lim\limits_{z\rightarrow0^+}m_F(z)=\dfrac{1}{1-c}\int\limits_{-\infty}^{+\infty}\dfrac{dH(\tau)}{\tau}\,,
\end{equation}
where the integral exists due to assumption (A3).

The function $m_F(z)$ is analytic in $\mathbbm{C}^+$ thus we take the derivative with respect to $z$ from both sides of equation (\ref{MP}) and get
\begin{equation}\label{derMP1}
m^\prime_F(z)=\int\limits_{-\infty}^{+\infty}\dfrac{\tau c(m_F(z)+zm^\prime_F(z))+1}{(\tau(1-c-czm_F(z))-z)^2}dH(\tau)\,.
\end{equation}
Rearranging terms in (\ref{derMP1}) we get
\begin{align}\label{derMP}
&m^\prime_F(z)\left(1-\int\limits_{-\infty}^{+\infty}\dfrac{cz\tau}{(\tau(1-c-czm_F(z))-z)^2}dH(\tau)\right)\\
&=\int\limits_{-\infty}^{+\infty}\dfrac{\tau cm_F(z)+1}{(\tau(1-c-czm_F(z))-z)^2}dH(\tau)\nonumber\,.
\end{align}
The right side of (\ref{derMP}) exists as $z\rightarrow0^+$ due to (\ref{modf}), $m_F(0)<\infty$, and the dominated convergence theorem, thus the left hand side must also exist as $z\rightarrow0^+$ and it implies that
\begin{equation}\label{limMP}
m^\prime_F(0)=\dfrac{1}{(1-c)^2}\left(cm_F(0)\int\limits_{-\infty}^{+\infty}\dfrac{dH(\tau)}{\tau}+\int\limits_{-\infty}^{+\infty}\dfrac{dH(\tau)}{\tau^2}\right)\,.
\end{equation}
The application of (\ref{lim0}) completes the proof of Theorem 3.1.\\

\textbf{Proof of Theorem 3.2}.
In order to prove Theorem 3.2 we need the following lemma.
\begin{lemma}\textbf{[Lemma B.26, Bai and Silverstein (2010)]}
 Let $\bA$ be a $p\times p$ nonrandom matrix and let $\bx=(x_1,\ldots,x_p)^\prime$ be a random vector with independent entries. Assume that $E(x_i)=0$, $E|x_i|^2=1$, and $E|x_i|^l\leq\nu_l$. Then, for any $k\geq1$,
 \begin{equation}\label{BSlem}
 E|\bx^\prime\bA\bx-\text{tr}(\bA)|^k\leq C_k\left((\nu_4\text{tr}(\bA\bA^\prime))^{\frac{k}{2}}+\nu_{2k}\text{tr}(\bA\bA^\prime)^{\frac{k}{2}}\right)\,,
 \end{equation}
 where $C_k$ is some constant which depends only on $k$.
 \end{lemma}
 Rubio and Mestre (2011) studied the asymptotics of the functionals $\text{tr}(\mathbf{\Theta}(\bS_n-z\bI)^{-1})$ for a deterministic matrix $\mathbf{\Theta}$ with bounded trace norm at infinity. It is noted that the results of Theorem 1 by Rubio and Mestre (2011) also hold under the assumption of the existence of $4$th moments which is weaker than the one given in the original paper.
 This statement is obtained by using Lemma B.26 of Bai and Silverstein (2010) on quadratic forms which we recall for presentation purposes as Lemma 7.1 above.

 In order to obtain the statement of Theorem 1 by Rubio and Mestre (2011) under the weaker assumption imposed on the moments, we replace Lemma 2 of Rubio and Mestre (2011) by Lemma 7.1 in the case of $k\geq1$. This implies that Lemma 3 of Rubio and Mestre (2011) holds also for $k\geq1$. Lemma 4 of Rubio and Mestre (2011) has already been proved under the assumption that there exist $4+\varepsilon$ moments. The last step is the application of Lemma 1, 2 and 3 of Rubio and Mestre (2011) with $k\geq1$. Finally, it can be easily checked that further steps of the proof of Theorem 1 by Rubio and Mestre (2011) hold under the existence of $4+\varepsilon$ moments.\\ A partial case of the result proved by Rubio and Mestre (2011) is summarized in Theorem 7.2.
\begin{theorem}{\textbf{[Rubio and Mestre (2011)]}}
Assume that (A2) and (A3) hold and additionally some nonrandom matrix $\mathbf{\Theta}$ has uniformly bounded trace norm at infinity then for $p/n\rightarrow c>0$ as $n\rightarrow\infty$
\begin{equation}\label{RM2011}
\left|\text{tr}\left(\mathbf{\Theta}(\bS_n-z\bI)^{-1}\right)-\text{tr}\left(\mathbf{\Theta}(x(z)\bSigma_n-z\bI)^{-1}\right)\right|\longrightarrow0~~\text{a. s.}\,
\end{equation}
where $x(z)$ is a unique solution in $\mathbbm{C}^+$ of the following equation
\begin{equation}\label{eq}
\dfrac{1-x(z)}{x(z)}=\dfrac{c}{p}\text{tr}\left(x(z)\bI-z\bSigma^{-1}_n\right)^{-1}\,.
\end{equation}
\end{theorem}

Next we prove Theorem 3.2 directly.
We investigate the asymptotic behavior of the following quantities
 \begin{align}
 &\gamma_1=\text{tr}(\bS^{-1}_n\bTheta)\label{g1},\\
 &\gamma_2=\dfrac{1}{p}||\bS_n^{-1}||^2_F=\dfrac{1}{p}\text{tr}(\bS^{-2}_n)\label{g2}\,.
 \end{align}
First, we consider the quantity $\gamma_1$ given in (\ref{g1}) and rewrite it as $\gamma_1(z)=\text{tr}((\bS_n-z\bI)^{-1}\bTheta)$ for all $z\in\mathbbm{C}^+$. Using Theorem 7.2 we get that for $\frac{p}{n}\rightarrow c\in(0, 1)$ as $n\rightarrow\infty$ holds
\begin{equation}\label{gamma1}
\left|\gamma_1(z)-\text{tr}(\bTheta(x(z)\bSigma_n-z\bI)^{-1})\right|\longrightarrow0~~\text{a. s.}\,
\end{equation}
and $x(z)$ is the unique solution in $\mathbbm{C}^+$ of the equation
\begin{equation}\label{eg1}
\dfrac{1-x(z)}{x(z)}=\dfrac{c}{p}\text{tr}(x(z)\bI-z\bSigma^{-1}_n)^{-1}\,.
\end{equation}
First, we show that
\begin{equation}\label{x0}
x(0)=\lim\limits_{z\rightarrow0^{+}}x(z)=1-c\,.
\end{equation}
If we assume that $x(0)=\infty$ then we get immediately the contradiction due to equation (\ref{eg1}). Similarly, from (\ref{eg1}) we conclude that $x(0)\neq0$. This implies that $0<x(0)<\infty$ and thus taking the limit $z\rightarrow0^+$ from both sides of (\ref{eg1}) we get (\ref{x0}).

Note that $\lim\limits_{z\rightarrow0^{+}}\gamma_1(z)=\gamma_1$ which together with (\ref{x0}) and (\ref{gamma1}) implies
\begin{equation}\label{gamma1}
\left|\gamma_1-\dfrac{1}{1-c}\text{tr}(\bSigma^{-1}_n\bTheta)\right|\longrightarrow0~~\text{a. s.}\,
\end{equation}
for $\dfrac{p}{n}\rightarrow c\in(0, 1)$ as $n\rightarrow\infty$.

Next, we prove the following statement
\begin{equation}\label{gamma2}
\Biggl|\gamma_2-\dfrac{1}{p}\dfrac{1}{(1-c)^2}\left(||\bSigma^{-1}_n||^2_{F}+\dfrac{c}{p(1-c)}(||\bSigma^{-1}_n||^2_{tr}\right)\Biggr|\longrightarrow0~~\text{a. s.}\,
\end{equation}
for $\frac{p}{n}\rightarrow c\in(0, 1)$ as $n\rightarrow\infty$.

Using the triangle inequality we rewrite the difference in (\ref{gamma2}) in the following way
\begin{align}\label{1step}
&\Biggl|\gamma_2-\dfrac{1}{p}\dfrac{1}{(1-c)^2}\left(||\bSigma^{-1}_n||^2_{F}+\dfrac{c}{p(1-c)}||\bSigma^{-1}_n||^2_{tr}\right)\Biggr|\leq\Biggl|\gamma_2-\psi\Biggr|\\
&+\Biggl|\psi-\dfrac{1}{p}\dfrac{1}{(1-c)^2}\left(||\bSigma^{-1}_n||^2_{F}+\dfrac{c}{p(1-c)}||\bSigma^{-1}_n||^2_{tr}\right)\Biggr|\nonumber\,,
\end{align}
where $\psi=\dfrac{1}{(1-c)^2}\int\limits_{-\infty}^{+\infty}\dfrac{1}{\tau^2}dH(\tau)+\dfrac{c}{(1-c)^3}\left(\int\limits_{-\infty}^{+\infty}\dfrac{1}{\tau} dH(\tau)\right)^2$ is given in Theorem 3.1. Next we show that the right side of (\ref{1step}) vanishes almost surely as $n\rightarrow\infty$. Using Theorem 3.1 we get
\begin{equation}\label{1term}
\Biggl|\gamma_2-\psi\Biggr|\longrightarrow0~\text{a. s.}~\text{for}~n\rightarrow\infty\,.
\end{equation}
Next, we show that the second nonrandom term in (\ref{1step}) approaches to zero as $n\rightarrow\infty$.
Using assumption (A1) it holds that $H_n(t)$ tends to $H(t)$ at all continuity points of $H(t)$. Thus,
\begin{equation}\label{tr1}
\dfrac{1}{p}||\bSigma^{-1}_n||^2_F=\dfrac{1}{p}\text{tr}(\bSigma^{-2}_n)=\dfrac{1}{p}\sum\limits_{i=1}^p\dfrac{1}{\tau^2_i}=\int\limits_{-\infty}^{+\infty}\dfrac{1}{\tau^2} dH_n(\tau)\stackrel{{n\rightarrow\infty}}{\longrightarrow}\int\limits_{-\infty}^{+\infty}\dfrac{1}{\tau^2} dH(\tau)\,,
\end{equation}
where the last integral in (\ref{tr1}) exists due to assumption (A3).
Similarly, it holds that
\begin{equation}\label{tr2}
\dfrac{1}{p^2}||\bSigma^{-1}_n||^2_{tr}\stackrel{{n\rightarrow\infty}}{\longrightarrow}\left(\int\limits_{-\infty}^{+\infty}\dfrac{1}{\tau} dH(\tau)\right)^2\,.
\end{equation}
Using (\ref{tr1}) and (\ref{tr2}) we get
\begin{equation}\label{2term}
\Biggl|\psi-\dfrac{1}{p}\dfrac{1}{(1-c)^2}\left(||\bSigma^{-1}_n||^2_{F}+\dfrac{c}{p(1-c)}||\bSigma^{-1}_n||^2_{tr}\right)\Biggr|\longrightarrow0~~\text{for}~n\rightarrow\infty
\end{equation}
As a result, (\ref{1term}) and (\ref{2term}) complete the proof of Theorem 3.2.\\

\textbf{Proof of Theorem 3.3}.
In the proof, we use a special case of the results derived by Rubio and Mestre (2011) which are summarized in Theorem 7.3.
\begin{theorem}{\textbf{[Rubio and Mestre (2011)]}}
Assume that (A2) and (A3) hold and additionally some nonrandom $n\times n$ matrix $\mathbf{\Theta}$ has uniformly bounded trace norm at infinity. Let $\bar{\bS}_n=\bx^\prime_n\bSigma_n\bx_n$. Then for $p/n\rightarrow c>0$ as $n\rightarrow\infty$
\begin{equation}\label{RM20112}
\left|\text{tr}\left(\mathbf{\Theta}(\bar{\bS}_n-z\bI)^{-1}\right)-\text{tr}\left(\mathbf{\Theta}\right)x(z)\right|\longrightarrow0~~\text{a. s.}\,
\end{equation}
where $x(z)$ is a unique solution in $\mathbbm{C}^+$ of the following equation
\begin{equation}\label{eq}
\dfrac{1+zx(z)}{x(z)}=\dfrac{c}{p}\text{tr}\left(x(z)\bI+\bSigma^{-1}_n\right)^{-1}\,.
\end{equation}
\end{theorem}
First we consider for $c>1$ the asymptotics of the quantity
\begin{equation}\label{k1}
\kappa_1=\dfrac{1}{p}\text{tr}(\bS^+_n)=\dfrac{c^{-1}}{n}\text{tr}(\bar{\bS}^{-1}_n)=\lim\limits_{z\rightarrow0^+}\dfrac{c^{-1}}{n}\text{tr}\left((\bar{\bS}_n-z\bI)^{-1}\right)\,,
\end{equation}
where the $n\times n$ matrix $\bar{\bS}_n$ is defines in Theorem 7.3. The second equality in (\ref{k1}) follows from the fact that the matrices $\bS_n$ and $\bar{\bS}_n$ possess the same nonzero eigenvalues (see, e.g., Silverstein (2009)). Note that the limit in (\ref{k1}) exists because $c>1$. Using Theorem 7.3 for $\bTheta=1/n\bI$ and setting $z\rightarrow0^+$ we get almost surely for $p/n\rightarrow c\in(1, +\infty)$
\begin{equation}\label{ask1}
\biggl|\kappa_1-c^{-1}x(0)\biggr|\rightarrow0~~\text{as}~n\rightarrow\infty\,
\end{equation}
and $x(0)=\lim\limits_{z\rightarrow0^+}x(z)$ satisfies the equation
\begin{equation}\label{eq2}
\dfrac{1}{x(0)}=\dfrac{c}{p}\text{tr}\left(x(0)\bI+\bSigma^{-1}_n\right)^{-1}\,.
\end{equation}
Note that $x(0)$ always exists in case of $c>1$. In order to find the asymptotics for the quantity $\kappa_1(\bTheta)=1/n\text{tr}(\bS^+_n\bTheta)$ we recall the properties of the generalized inverse. Namely, if $\bTheta$ is positive definite then $\bTheta^{1/2}\bS^{+}_n\bTheta^{1/2}$ is a generalized inverse of $\bTheta^{-1/2}\bS_n\bTheta^{-1/2}$ which has the same non-zero eigenvalues as $\bar{\bS}^*=\bx^\prime\bSigma_n^{1/2}\bTheta^{-1} \bSigma_n^{1/2}\bx$. Consequently, since $tr(\bS^{+}\bTheta)=tr(\bTheta^{1/2}\bS^{+}_n\bTheta^{1/2})$ we get the same result as in case of $\kappa_1=1/n\text{tr}(\bS^{+}_n)$ with $\bSigma_n^{-1}$ be replaced by
\begin{equation}
(\bSigma_n^{1/2}\bTheta^{-1}\bSigma_n^{1/2})^{-1}=\bSigma_n^{-1/2}\bTheta\bSigma_n^{-1/2}\,.
\end{equation}

The second quantity of interest is
\begin{equation}\label{k2}
\kappa_2=\dfrac{1}{p}||\bS^+_n||_F^2=\dfrac{c^{-1}}{n}||\bar{\bS}^{-1}_n||_F^2=c^{-1}\lim\limits_{z\rightarrow0^+}\dfrac{\partial}{\partial z}\dfrac{1}{n}\text{tr}\left((\bar{\bS}_n-z\bI)^{-1}\right)\,.
\end{equation}
Again, we use Theorem 7.3 for $\dfrac{1}{n}\text{tr}\left((\bar{\bS}_n-z\bI)^{-1}\right)$, calculate the derivative with respect to $z$ and set it to zero. As a result, we obtain the following identity for $p/n\rightarrow c\in(1,+\infty)$
\begin{equation}\label{ask2}
\biggl|\kappa_2-c^{-1}x^\prime(0)\biggr|\rightarrow0~~\text{as}~n\rightarrow\infty\,,
\end{equation}
 where
\begin{equation}\label{xprime0k2}
\dfrac{1}{x^\prime(0)}=\dfrac{1}{x^2(0)}-\dfrac{c}{p}\text{tr}\left(x(0)\bI+\bSigma^{-1}_n\right)^{-2}\,
\end{equation}
 and $x(0)$ satisfies the equation (\ref{eq2}). From (\ref{ask1}) and (\ref{ask2}) follows the statement of Theorem 3.3.\\

\textbf{Proof of Proposition 3.1}.
Using Theorem 3.3 we get that $\dfrac{1}{p}\boldsymbol{\eta}^\prime\bS^{+}_n\boldsymbol{\xi}=
\dfrac{1}{p}\text{tr}(\boldsymbol{\xi}\boldsymbol{\eta}^\prime\bS^{+}_n)\underset{\text{a.s.}}{\rightarrow}c^{-1}y(\bxi\boldsymbol{\eta}^\prime)$, where $y(\bxi\boldsymbol{\eta}^\prime)$ satisfies the equation
\begin{equation}\label{xixi}
\dfrac{1}{y(\boldsymbol{\xi}\boldsymbol{\eta}^\prime)}=\dfrac{c}{p}\text{tr}\left[\left(\bSigma_n^{-1/2}
\boldsymbol{\xi}\boldsymbol{\eta}^\prime\bSigma_n^{-1/2}+y(\boldsymbol{\xi}\boldsymbol{\eta}^\prime)\bI\right)^{-1}\right]\,.
\end{equation}

Using the Sherman-Morrison formula (see, e.g., Horn and Johnson (1986)) we can rewrite the the right-hand side of (\ref{xixi}) in the following way
\begin{equation}\label{xixi2}
\dfrac{c}{p}\text{tr}\left[\dfrac{1}{y(\boldsymbol{\xi}\boldsymbol{\eta}^\prime)}\bI-\dfrac{1}{y^2(\boldsymbol{\xi}\boldsymbol{\eta}^\prime)}\dfrac{\bSigma_n^{-1/2}
\boldsymbol{\xi}\boldsymbol{\eta}^\prime\bSigma_n^{-1/2}}{1+\dfrac{1}{y(\boldsymbol{\xi}\boldsymbol{\eta}^\prime)}
\boldsymbol{\eta}^\prime\bSigma_n^{-1}\boldsymbol{\xi}}\right]
=\dfrac{c}{y(\boldsymbol{\xi}\boldsymbol{\eta}^\prime)}\left(1-\dfrac{
\boldsymbol{\eta}^\prime\bSigma_n^{-1}\boldsymbol{\xi}}{y(\boldsymbol{\xi}\boldsymbol{\eta}^\prime)+
\boldsymbol{\eta}^\prime\bSigma_n^{-1}\boldsymbol{\xi}}\right)\,.
\end{equation}

Combining (\ref{xixi}) with (\ref{xixi2}) and multiplying both sides by $y(\boldsymbol{\xi}\boldsymbol{\eta}^\prime)$, we get
\begin{equation}\label{xixiend}
1-c^{-1}=\dfrac{\boldsymbol{\eta}^\prime\bSigma_n^{-1}\boldsymbol{\xi}}{y(\boldsymbol{\xi}\boldsymbol{\eta}^\prime)+
\boldsymbol{\eta}^\prime\bSigma_n^{-1}\boldsymbol{\xi}}\,,
\end{equation}
which is a linear equation in $y(\boldsymbol{\xi}\boldsymbol{\eta}^\prime)$ with the solution given by
\begin{equation}\label{xixisol}
y(\boldsymbol{\xi}\boldsymbol{\eta}^\prime)=\dfrac{c^{-1}}{1-c^{-1}}\boldsymbol{\eta}^\prime\bSigma^{-1}_n\boldsymbol{\xi}
=\dfrac{1}{c-1}\boldsymbol{\eta}^\prime\bSigma^{-1}_n\boldsymbol{\xi}\,.
\end{equation}
The last equality finishes the proof of Proposition 3.1.\\

\textbf{Proof of Proposition 4.1}.
We write the proof in the case $c<1$ because the case $c>1$ is already handled in Section 3.2.
\begin{enumerate}[label=\roman{*}).]
\item Consider the functionals $\widehat{\alpha}^*_n(\sigma\boldsymbol{\Pi}_0)$ and $\widehat{\beta}^*_n(\sigma\boldsymbol{\Pi}_0)$. Using (\ref{ea}) and (\ref{eb}) it holds that
\begin{align}
\label{aaa}
\widehat{\alpha}^*_n(\sigma\boldsymbol{\Pi}_0)&=\widehat{\alpha}^*_n(\boldsymbol{\Pi}_0)\\
\label{bbb}
\widehat{\beta}^*_n(\sigma\boldsymbol{\Pi}_0)&=\dfrac{1}{\sigma}\widehat{\beta}^*_n(\boldsymbol{\Pi}_0)\,.
\end{align}
Putting (\ref{aaa}) and (\ref{bbb}) together with $\sigma\boldsymbol{\Pi}_0$ in (\ref{olse}) completes the proof of the first part of Proposition 4.1.\\
\item
From Corollary 3.1 it holds that
\begin{equation}\label{ainftys}
\left|\widehat{\alpha}^*_n(1/p\bI)-\ta^*(1/p\bI)\right|\longrightarrow0~~\text{a. s.}~~\text{for}~n\rightarrow\infty\,
\end{equation}
with
\begin{equation}\label{a_ass}
\ta^*(1/p\bI)=(1-c)\dfrac{||\bSigma_n^{-1}||^2_F||1/p\bI||^2_F-\left(\text{tr}(\bSigma^{-1}_n1/p\bI)\right)^2}{\left(||\bSigma_n^{-1}||^2_F+\dfrac{c}{p(1-c)}||\bSigma_n^{-1}||_{tr}^2\right)||1/p\bI||^2_F-\left(\text{tr}(\bSigma^{-1}_n1/p\bI)\right)^2}
\end{equation}
Using (\ref{a_ass}) and noting that $\bSigma^{-1}_n=\sigma\bI$ we get that
\begin{equation}\label{a0}
\ta^*(1/p\bI)=0.
\end{equation}
Similarly,
\begin{equation}\label{binftys}
\left|\widehat{\beta}^*_n(1/p\bI)-\tb^*(1/p\bI)\right|\longrightarrow0~~\text{a. s.}~~\text{for}~n\rightarrow\infty\,
\end{equation}
and
\begin{equation}\label{b_ass}
\tb^*(1/p\bI)=\dfrac{\text{tr}(\bSigma^{-1}_n1/p\bI)}{||1/p\bI||^2_F}\left(1-\dfrac{\alpha^*(1/p\bI)}{1-c}\right)\,,
\end{equation}
which together with (\ref{a0}) and $\bSigma^{-1}_n=\sigma\bI$ implies that
\begin{equation}\label{b0}
\tb^*(1/p\bI)=p\sigma.
\end{equation}
The equalities (\ref{a0}) and (\ref{b0}) with (\ref{olse}) complete the proof of the second part of Proposition 4.1.\\
\item
From (\ref{a_as}) and (\ref{b_as}) in Corollary 3.1 it follows that

\begin{equation}\label{a11}
\alpha^*_n(\bSigma_n^{-1})=0~~\text{and}~~\beta^*_n(\bSigma_n^{-1})=1\,
\end{equation}
and, hence,
\begin{equation}\label{a00}
\widehat{\alpha}^*_n(\bSigma_n^{-1})\longrightarrow0~~\text{and}~~\widehat{\beta}^*_n(\bSigma_n^{-1})\longrightarrow1~\text{a.s.}~\text{for}~n\rightarrow\infty\,.
\end{equation}
 Substituting (\ref{a00}) in (\ref{olse}) completes the proof of the third part of Proposition 4.1.
\end{enumerate}

\end{document}